\newif\iffinal
\finalfalse	% Not final version
\finaltrue	% Final version

\documentclass[letterpaper,12pt,reqno]{amsart}
\RequirePackage[utf8]{inputenc}
\usepackage[portrait,margin=2.6cm]{geometry}
\iffinal\else\usepackage[notref,notcite]{showkeys}\fi
\usepackage{mathrsfs,soul}
\usepackage{hyperref}
\usepackage[foot]{amsaddr}
\usepackage{amssymb,amsthm,amsfonts,amsbsy,latexsym,dsfont}
\usepackage[textsize=small]{todonotes}       % includes TO DO LIST
\usepackage{graphicx}
\usepackage[numeric,initials,nobysame]{amsrefs}
\usepackage{upref,setspace}
\usepackage{xcolor,colortbl}
\usepackage{enumerate}
\newenvironment{enumeratei}{\begin{enumerate}[\upshape (i)]}{\end{enumerate}}
\newenvironment{enumeratea}{\begin{enumerate}[\upshape (a)]}{\end{enumerate}}

\usepackage{paralist}

\numberwithin{equation}{section}
\numberwithin{figure}{section}
\numberwithin{table}{section}

\newtheorem{thm}{Theorem}[section]
\newtheorem{lem}[thm]{Lemma}

\newtheorem{cor}[thm]{Corollary}

\newtheorem{prop}[thm]{Proposition}

\newtheorem{lemma}[thm]{Lemma}

\theoremstyle{definition}
\newtheorem{rem}{Remark}

\newcommand{\eps}{\varepsilon}

\newcommand{\set}[1]{\left\{#1\right\}}

\newcommand{\ie}{\emph{i.e. }}

\newcommand{\eg}{\emph{e.g. }}

\newcommand{\equald}{\stackrel{\mathrm{d}}{=}}

\newcommand{\convas}{\stackrel{\mathrm{a.s.}}{\longrightarrow}}

\def\qed{ \hfill $\blacksquare$}
\newcommand{\cA}{\mathcal{A}}
\newcommand{\cF}{\mathcal{F}}
\newcommand{\cI}{\mathcal{I}}

\newcommand{\cT}{\mathcal{T}}\newcommand{\cU}{\mathcal{U}}

\newcommand{\bE}{\mathbb{E}}

\newcommand{\bL}{\mathbb{L}}

\newcommand{\bR}{\mathbb{R}}

\newcommand{\bZ}{\mathbb{Z}}
\newcommand{\rE}{\mathrm{E}}

\DeclareMathOperator{\E}{\mathbb{E}}
\DeclareMathOperator{\pr}{\mathbb{P}}

 \DeclareMathOperator{\BP}{pdBP}

\newcommand{\convd}{\stackrel{d}{\longrightarrow}}
\newcommand{\convp}{\stackrel{P}{\longrightarrow}}

\usepackage{fourier}
\usepackage{pdfsync}

\begin{document}

\title[Community modulated recursive trees]{Community modulated recursive trees and population dependent branching processes}

\date{}
% Secondary: ;}

\author[Bhamidi]{Shankar Bhamidi}
\author[Fan]{Ruituo Fan}
\author[Fraiman]{Nicolas Fraiman}
\author[Nobel]{Andrew Nobel}

\address{Department of Statistics and Operations Research, 304 Hanes Hall, University of North Carolina, Chapel Hill, NC 27599}
\email{bhamidi@email.unc.edu,badhron@email.unc.edu,fraiman@email.unc.edu,nobel@email.unc.edu}
\maketitle
\begin{abstract}
	We consider random recursive trees that are grown via community modulated schemes that involve random attachment or degree based attachment. The aim of this paper is to derive general techniques based on continuous time embedding to study such models. The associated continuous time embeddings are {\bf not} branching processes: individual reproductive rates at each time $t$ depend on the composition of the entire population at that time. Using stochastic analytic techniques we show that various key macroscopic statistics of the continuous time embedding stabilize, allowing asymptotics for a host of functionals of the original models to be derived.  
\end{abstract}

\section{Introduction} 

The subject of this paper is the analysis of random recursive trees that are grown via community modulated 
schemes involving random attachment or degree based attachment.
This topic lies at the intersection of two areas of modern probabilistic combinatorics: random recursive trees and community detection for networks. 
Before describing the models of primary interest, we give some preliminaries. 

A rooted tree on $n$ vertices labeled $\{ 1,2,...,n \}$ is called recursive if $1$ is the root, and for each $2 \leq i \leq n$, the (unique) path from the root to vertex $i$ has increasing labels. The name ``recursive'' suggests that these trees can be constructed recursively, by adding a new vertex at each time step, and the labels can be viewed as birth orders. Recursive trees have been studied for decades (\eg see Smythe and Mahmoud \cite{smythe1995survey} for an early survey), with applications to epidemics \cite{moon1974distance}, pyramid schemes \cite{gastwirth1977probability}\cite{gastwirth1984two}, convex hull algorithms \cite{mahmoud1992asymptotic}, and modeling family trees of preserved copies of ancient or medieval texts \cite{najock1982number}, where vertices represents people or texts that arrive chronologically and are labeled by time. 

In many studies of rooted trees, the uniform recursive tree (URT) model is used, \ie, a tree is chosen uniformly at random from all recursive trees of a fixed size. From a statistical point of view, URTs serve are a natural null model for recursive trees.  Recently, a wide array of extensions of URTs have been considered. For example, preferential attachment models \cite{barabasi1999emergence}, which give rise to power law degree distributions, are favored in the complex network community where real data exhibits heavy-tailed degrees.  A recent multitype extension of preferential attachment based on the genealogy of a multitype branching process has been considered in \cite{rosengren2017multi}. There are also variants of URTs that introduce choices to the attachment rule \cite{d2007power}\cite{mahmoud2010power}. Instead of choosing one existing vertex, these models choose $k$ previous vertices (with or without replacement) as candidates, and connect the new vertex to one of them based on certain optimization criterion. Related to these models is the so called scaled attachment random recursive tree (SARRT) \cite{devroye2012depth}, where at the $n$-th step the new vertex $n$ is connected to vertex $\lfloor n X_n \rfloor$, with $X_1,X_2,...$ being a sequence of i.i.d. random variables taking value in $[0,1)$. When the $X_i$ follow a uniform distribution, SARRT reduces to the usual URT. 

However, all of these alternatives are still homogeneous in the sense that attachments are made based on the same rule, no matter how complicated. In this paper, we introduce an alternative to URT that models heterogeneity in the spirit of the well-studied \emph{stochastic block model}. Recall that in the stochastic block model we start with $n$ vertices that are partitioned into two or more classes (often referred to as ``communities''), and that each pair of vertices are connected independently at random with probability depending only on their class membership.  Just as the stochastic block model can be seen as a heterogeneous alternative to the Erd\H{o}s-R\'{e}nyi model (where each pair of vertices are connected independently with the same probability), our model, which we shall refer to as Community Modulated Recursive Tree (CMRT), generalizes URT in a similar way to allow for latent class labels.

\subsection{Model formulation}
Here we describe the CMRT canonical model. Extensions of this model to more than two types as well as other schemes of attachment including preferential attachment are explored in Section \ref{sec:exten}. First we fix some terminology for rooted trees. Since they are connected (undirected) graphs with no cycles, there exists a unique path from any given vertex $v$ to the root. Any vertex $u \neq v$ in this path is called an antecedent of $v$, and $v$ its descendant; if $u$ and $v$ are adjacent (\ie there exists an edge between them), we call $u$ the parent of $v$ and $v$ its child. We also say that $v$ is connected to $u$ if they are adjacent. For recursive trees with multiple vertices identified as roots, these terms are defined with respect to the nearest root. 

Recall that a uniform recursive tree (URT) can be constructed as follows: starting with a root vertex labeled 1, at each time $n$ choose an existing vertex uniformly at random and add a new vertex labeled $n$ connected to this chosen vertex. Note that this definition yields a growing tree-valued process, which we shall denote by $\set{\cU_n}_{n \geq 1}$, where $\cU_n$ is the random recursive tree given by this process when it reaches size $n$. We shall refer to $\set{\cU_n}_{n\geq 1}$ as a URT and $\cU_n$ as a URT of size $n$.

Now we are ready to introduce Community Modulated Recursive Tree (CMRT). This is the ``canonical'' model in the sense that it is easiest to describe both in discrete time and through a continuous time embedding that allows for its analysis.  Once the reader is familiar with this model, understanding extensions becomes easy.  Throughout this paper we shall always refer to classes of vertices in random trees as ``types''. For simplicity we first state results for a two-type model, and postpone the general case to Section \ref{sec:exten}. Specifically, assume that all vertices are of either type $A$ or type $B$. 
To construct a Community Modulated Recursive Tree (CMRT) we begin with a single vertex 1 of
type $A$ and a single vertex 2 of type $B$ that are connected by an edge.  This initial tree is the CMRT of 
size 2.  Vertices 1 and 2 as roots of type $A$ and type $B$, respectively.  	At each subsequent stage $n \geq 3$, 
a new vertex labeled $n$ is added to the given tree of size $n-1$ in a three-step process:
\begin{enumeratea}
\vskip.1in
\item[(1)] Vertex $n$ flips a $p$-coin to determine if it is of type $A$ (with probability $p \in (0,1]$) 
or of type $B$ (with probability $1-p$). 
\vskip.1in
\item[(2)] Vertex $n$ flips a $q$-coin to decide if it will connect to vertices of its own type (with probability $q \in [0,1]$) 
or to vertices of the other type (with probability $1-q$). 
	\vskip.1in
\item[(3)] Vertex $n$ chooses an existing vertex of the type selected in step (2) uniformly at random  
and connects to it, forming a CMRT of size $n$.
\end{enumeratea} 
Note that in (a) we do not lose any generality by excluding the case where $p=0$, as one can switch the type labels
if necessary.  Adopting the terminology of stochastic block models, we call the CMRT assortative when 
$q \geq 1/2$ and disassortative when $q < 1/2$. In general, the probability that a vertex chooses to connect to
a vertex of the same type may vary across types.  We assume here that they are the same for simplicity; 
the general case is discussed in Section \ref{sec:exten}.

We refer to the random tree $\cT_n$ produced at the $n$-th stage of the procedure above as a Community 
Modulated Recursive Tree (CMRT) of size $n$, and to the sequence $\set{\cT_n}_{n\geq 2}$ as a CMRT.

\begin{figure}
	\centering
	\includegraphics[width=0.7\textwidth]{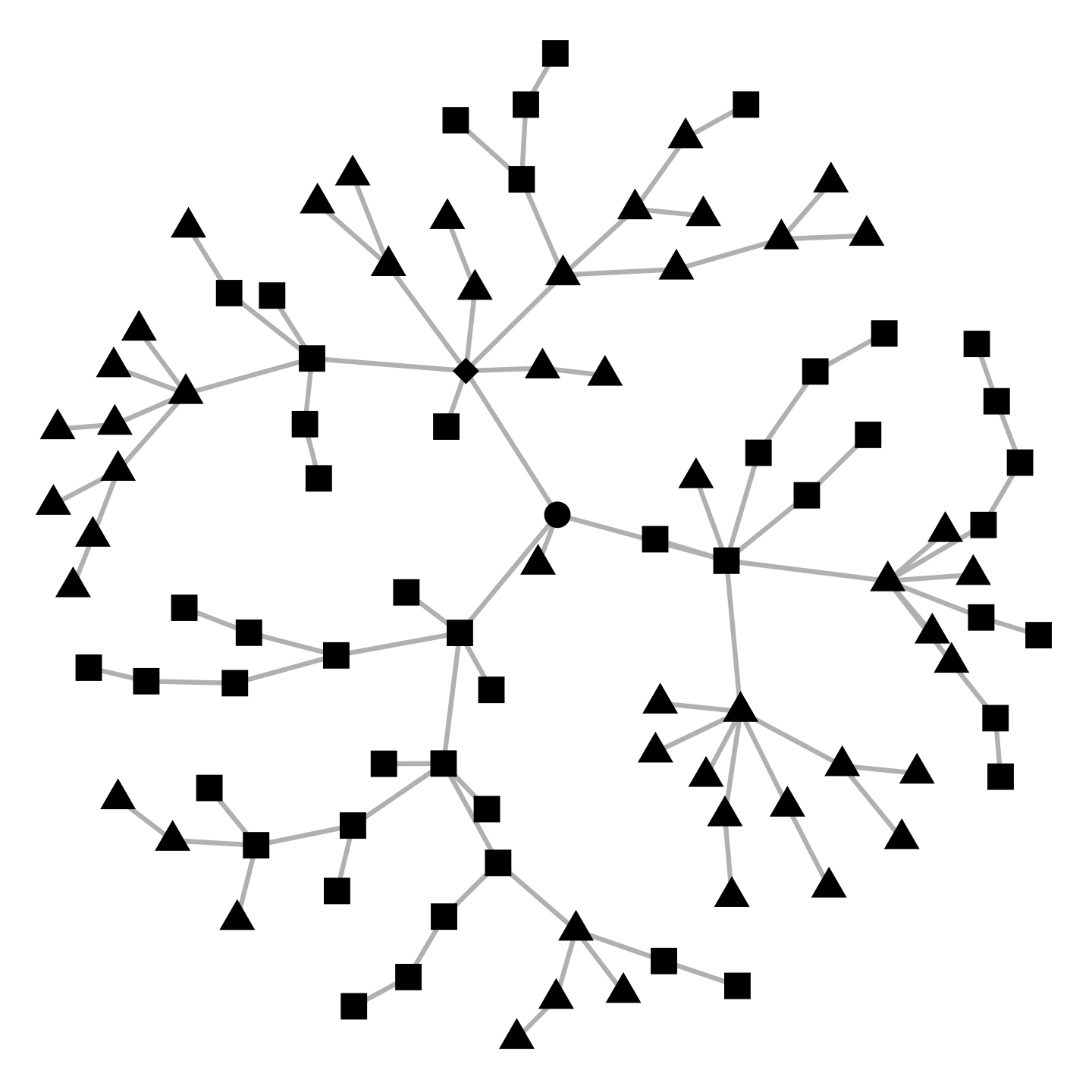}
	\caption{A (two-type) Community Modulated Recursive Tree of size 100, with parameters $p=0.5$ and $q=0.8$. The type $A$ and $B$ roots are respectively plotted as circle and diamond, and other type $A$ and $B$ vertices as squares and triangles.}
\end{figure}

\subsection{Special cases}
\label{sec:special}
\begin{enumeratei}
\vskip.1in
	\item When $p=1$ and $q=1$, the CMRT becomes a URT consisting of type $A$ vertices plus an additional type $B$ vertex (\ie the type $B$ root) connecting to the type $A$ root. Thus one may consider this case to be degenerate. 
	\vskip.1in
	\item 
	\label{item:spec-construct}
	When $p=1$ and $q\neq 1$, all new vertices are of type $A$ and have a fixed probability $1-q$ to connect to the unique type $B$ vertex.  This case too is degenerate in the sense that one can construct such a CMRT from a URT. To do so, take a URT consisting only of type $A$ vertices, and add a single type $B$ vertex connected to the root of the URT.  Treat this type $B$ vertex as the type $B$ root. Then remove each edge between type $A$ vertices independently with probability $1-q$, and for each edge removed in this way, add an edge between the corresponding child and the type $B$ root. It is easy to see that the recursive tree formed this way has the same distribution as a CMRT with $p=1$ and $q \neq 1$. This construction will be useful in some of the proofs.
	\vskip.1in
	\item When $q=1$, the CMRT can be partitioned into two disjoint, rooted subtrees associated with the type $A$ and type $B$ roots.  Furthermore, each subtree, conditioned on its size, has the same distribution as a URT. 
However, in the absence of vertex labels, identifying these subtrees from an observed CMRT of finite size is non-trivial.
Although not fully comparable, we note by way of comparison that in a stochastic block recovering the block assignment for each vertex is trivial if each block is connected and the probability to connect to vertices in other blocks is zero.
	\vskip.1in
	\item When $q=0$, each new vertex connects to a vertex of the opposite type. Thus for any path in the graph, adjacent vertices will have alternating types.
\end{enumeratei}

\subsection{Extensions}
\label{sec:exten}
Here we describe several natural extensions of the canonical model above.

\subsubsection{General CMRT}
One may readily extend CMRTs to $K>2$ types. As in the $K =2$ case we begin with one vertex of each type 
to which subsequent vertices may connect.  The initial tree $n = K$ is a URT of size $K$, the vertices of which
are assigned to the $K$ types based on uniform random permutation. 
We shall refer to these vertices as the roots of each type. 
At each subsequent stage $n \geq 3$, 
a new vertex labeled $n$ is added to the given tree of size $n-1$ in a three-step process:
\begin{enumeratea}
\vskip.1in
\item[(1)]
Vertex $n$ is assigned to type $i \in \{ 1, \ldots, K \}$ with probability $p_i$, where we assume without
loss of generality that $p_1 = \max_{1 \leq i \leq K} p_i > 0$.
\vskip.1in
\item[(2)] 
Once its type has been determined, vertex $n$ chooses to connect to a vertex of type $j \in \{ 1, \ldots, K \}$ 
with probability $q_{ij}$. 
\vskip.1in
\item[(3)] 
Vertex $n$ chooses an existing vertex of the type selected in step (2) uniformly at random  
and connects to it, forming a CMRT of size $n$.
\end{enumeratea} 
We refer to the sequence $\set{\cT_n}_{n \geq K}$ of stochastic trees produced in this way as 
a Community Modulated Recursive Tree (CMRT) with $K$ types.

\subsubsection{Community Weighted Recursive Tree (CWRT)}
There are other variants of multitype recursive trees that could be of potential interest.  
For example, one can generate CMRT with $K$ types following steps (1) and (2) above, but let the attachment probabilities in step (3) depend on weights for transitions between types.  
In more detail, let $\{ \omega_{ij} : 1 \leq i,j \leq K \}$ be positive weights, and let the incoming vertex $n$ of type
$i$ attach to an existing vertex of type $j$ with probability
\[
\frac{\omega_{ij}}{\sum_{l=1}^K n_l \omega_{il}},
\]
where $n_l$ is the number of existing type $l$ vertices in the tree. 
We refer to the sequence $\set{\cT_n}_{n \geq K}$ of stochastic trees produced in this way as 
a Community Weighted Recursive Tree (CWRT) with $K$ types.

\subsubsection{Community Modulated Preferential Attachment (CMPA)}
We can also define a ``preferential attachment'' version of CMRT, which we shall refer to as Community Modulated Preferential Attachment (CMPA) model. To construct a CMPA tree with $K$ types, we follow the same rule as that of a CMRT with $K$ types, except for the attachment step.  Suppose that a vertex $n$ of type $i$ that has chosen to connect to a vertex of type $j$.  Rather than choosing uniformly among vertices of type $j$, vertex $n$ connects to a type $j$ vertex $v$ with probability proportional to $D_v + \alpha_{ij}$, where $D_v$ is the out-degree of vertex $v$, and $\alpha_{ij}$ are positive numbers. 
We refer to the sequence $\set{\cT^*_n}_{n \geq K}$ as a Community Modulated Preferential Attachment (CMPA) tree with $K$ types.

\subsection{Contributions of the paper} 
In addition to the formulation of the models described above, the main aim of this paper is to develop a common set of tools based on continuous time embedding of the associated discrete time processes that allows for asymptotics in the large network $n \to \infty$ limit to be derived. 
The starting point, described in more detail in Section \ref{sec:ct-em}, is the embedding of recursive trees into a continuous time process.  We emphasize that the continuous time process is {\bf not} a branching process, in the sense that the instantaneous offspring growth rates of each vertex are intricately tied to the frequency of vertex types at that time. Using techniques from stochastic analysis, one can show that key macroscopic properties of the continuous time embedding (properly normalized) stabilize, thus allowing one to read off asymptotics for the discrete time processes. 
In particular, we shall study the limiting degree distribution, depth and height of these recursive trees, and establish the corresponding asymptotics.

\subsection{Organization of the paper} 
The next section describes the population dependent branching process relevant for the canonical model in Section \ref{sec:ct-em}.  
Analogous extensions of the embedding for other models including the community modulated preferential attachment are described in the proofs (Section \ref{sec:proof}).  Section \ref{sec:res} contains 
the statement of our main theoretical results. In Section \ref{sec:related} we provide a brief discussion on 
the relevance of this paper as well as related work. Section \ref{sec:proof} contains proofs of all the results.

\section{Continuous Time Embedding}
\label{sec:ct-em}
In this section we describe the continuous time embedding of recursive trees, which is the major tool used in our 
theoretical analyses. 
To fix ideas we describe an embedding for the canonical two-type Community Modulated Recursive Tree introduced above.  
In order to state the basic embedding result, we introduce a (two-type) continuous time process that 
we refer to as a {\bf population dependent} branching process ({\bf pdBP}). 
We refer to the vertices in the continuous time process as ``individuals'' to differentiate them from vertices in the corresponding 
discrete time random tree. 
This process is not a (homogeneous or inhomogeneous) branching process in the usual sense as the reproduction processes 
of different individuals are not independent. 
However, lengths of the time intervals between consecutive births are still exponentially distributed, and the parameters (expected number of births in a unit time interval) of these exponential variables are referred to as ``rates'', similar to that of the usual branching process.

Each individual in the population dependent branching process has type $A$ or $B$ (not both) and lives forever. 
The process is initialized at time $t=0$ with two individuals, one of each type, 
that we will refer to as ``ancestors''. 
At times $t > 0$ each existing individual in the process
gives birth to new individuals (offspring) of type $A$ or $B$ at rates that are specified below.

For each $t \geq 0$ let $n_A(t)$ and $n_B(t)$ denote the number of individuals of type $A$ and $B$, respectively, in the process, 
initialized with $n_A(0) = n_B(0) =1$.  Reproduction rates for individuals are as follows.  For each $t \geq 0$

\vskip.1in

\begin{enumeratei}

\item 
Each type $A$ individual gives birth to type $A$ individuals at rate $r_{AA}(t) = q$, and to type $B$ individuals at rate $r_{AB}(t) = (1-p)(1-q)/p$. 

\vskip.1in

\item
Each type $B$ individual gives birth to type $A$ individuals at rate
		\[
		r_{BA}(t) = \frac{n_A(t)}{n_B(t)} \cdot (1-q)
		\]
and to type $B$ individuals at rate
		\[
		r_{BB}(t) = \frac{n_A(t)}{n_B(t)} \cdot \frac{ q(1-p) }{p}.
		\]

\end{enumeratei}   

\vskip.1in

We introduce some notation that will be useful in what follows. 
For $t \geq 0$, let $n(t):=n_A(t) + n_B(t)$ be the total number of individuals alive at time $t$ and $\cF(t)$ be the $\sigma$-field generated by the process until time $t$. 
Also, denote by $\set{\cF(t):t \geq 0}$ the natural filtration of the process. 
For each $t \geq 0$ the genealogical structure of the population dependent branching process
can be described by a recursive tree in which each individual corresponds to a vertex.  In detail, 
the $A$ and $B$ ancestors are labeled 1 and 2, respectively, and are connected by an edge.
The vertex for an individual is connected by edges to the vertices of its offspring, and is labeled
by the absolute birth order of the individual in the overall process.
Let $\BP(t)$ denote the recursive tree capturing the genealogical structure of the process until 
and including time $t$.  The next lemma establishes a close connection between the continuous
time process $\set{\BP(t) : t \geq 0}$ and the CMRT.  The proof of the lemma is given in Section \ref{sec:proof}.

\begin{lemma}\label{lem:embed}
Let $\set{\BP(t) : t \geq 0}$ be the continuous time process described above.  For $n \geq 2$ define the stopping time
$T_n = \inf\set{ t \geq 0 : n(t) =n }$. 
	Then $\BP(T_n)  \stackrel{d}{=} \cT_n$ where $\set{\cT_n}_{n\geq 2}$ is a CMRT. 
	In fact,  
	$\set{\BP(T_n)}_{n\geq 2}  \stackrel{d}{=} \set{\cT_n}_{n\geq 2}$ as processes.  
\end{lemma}

\section{Results}
\label{sec:res}

\subsection{Asymptotics for CMRT}
\label{sec:cmrt-res}
Our first result concerns the limiting distribution of the out-degrees, \ie, the number of children, of vertices in the canonical CMRT. 
If one regards each edge in a rooted tree as directed from parent to child, then this definition is consistent with that of out-degree 
in directed graph.

\begin{thm}\label{thm:degree-dist}
Let $\{ \cT_n : n \geq 1 \}$ be the canonical CMRT, and for each $k \geq 0$, 
let $N_k(n)$ denote the number of vertices with out-degree $k$ in $\cT_n$. 
Then for each fixed $k \geq 0$ the ratio $N_k(n) / n$ converges in probability to $p_k$ where
	\begin{equation}
	\label{eqn:pk-def}
	p_k := \frac{p}{ 1 + r^*_A } \left( \frac{r^*_A}{ 1 + r^*_A } \right) ^k + 
	\frac{ 1-p }{ 1 + r^*_B } \left( \frac{r^*_B}{ 1 + r^*_B } \right) ^k
	\end{equation}
	and
	\[
	r^*_A:= q + \frac{ (1-p)(1-q) }{p} 
	\quad \text{and} \quad 
	r^*_B:= \frac{p (1-q) }{1-p} + q.
	\]
	When $p=1$ the second term in $p_k$ should be interpreted as 0.
\end{thm}

\begin{rem}
	The limiting degree distribution is a mixture of two geometric distribution shifted to the left by 1, 
	with parameters $1/(1 + r^*_A)$ and $1/(1 + r^*_B)$.  This distribution is identical to the limiting
	degree distribution of a URT \cite{moon1974distance}
	if and only if $p=1/2$ or $q=1$.
	Moreover, when $p \neq 1/2$ and $q \neq 1$, the limiting proportion of leaves $p_0$ satisfies
	\[
	p_0 - \frac{1}{2} = 
	\frac{p}{ 1 + r^*_A } + \frac{ 1-p }{ 1 + r^*_B } -\frac{1}{2} =
	\frac{ (2p-1)^2 (1-q)^2 }{ 2(1-(1-2p)^2 q^2) } > 0.
	\]
	Thus in the case where the limiting degree distribution of the CMRT differs from that of the URT, the former has a larger proportion of leaves. 
	In fact, $p_0 \to 1$ as $p \to 1$ and $q \to 0$.  In this extreme case, $\cT_n$ looks like a ``star''  in which $n-1$ type $A$ vertices 
	are connect to the type $B$ root.
\end{rem}

To better understand Theorem \ref{thm:degree-dist}, we outline the following heuristic calculations based on the continuous time embedding. 
Consider the case where $p<1$.  Note that the proportions of type $A$ and $B$ vertices in the CMRT $\set{\cT_n}_{n \geq 2}$
converge almost surely to $p$ and $1-p$ respectively as $n \to \infty$ by the strong law of large numbers. 
If one regards $n_A(t)$ and $n_B(t)$ as deterministic differentiable functions, formal solution of the rate equations (i) and (ii) in the
definition of the process $\BP(t)$ yields $n_A(t) = e^t$ and $n_B(t) = (1-p) e^t / p$.  This suggests that, in analogy with the CMRT,
\[
\frac{n_A(t)}{n_B(t)} \approx \frac{p}{1-p}
\]
when $t$ is large.  One may make this rigorous using stopping time arguments. 
Plugging this approximation into (ii) yields corresponding approximations for the transition rates of type $B$ individuals
\[
r_{BA}(t) \approx \frac{p}{1-p}\cdot (1-q),
\qquad 
r_{BB}(t) \approx \frac{p}{1-p}\cdot \frac{ q(1-p) }{p}=q.
\]
Note that these approximations, and the transition rates (i) for type $A$ individuals are independent of $t$.
Combining the transition rates for type $A$ we get $r_{AA}(t)+r_{AB}(t) = r^*_A$.  Similarly, combining the
approximate rates for type $B$ individuals yields $r_{BA}(t)+r_{BB}(t) \approx r^*_B$.

These calculations suggest that one can approximate the population-dependent branching process by a time 
homogeneous multi-type branching process with types $A, B$, and rates $r^*_A, r^*_B$. 
Results of Jagers and Nerman \cite{jagers1996asymptotic}, applied to the approximation, show that 
the type of a randomly chosen individual $v$ and its age have joint limiting distribution
$\pi \times \text{Exponential}(1)$, where  
$\pi$ is a discrete measure on $\set{A,B}$ with $\pi(A)=p$ and $\pi(B)=1-p$ and $\times$ denotes product measure. 
As a type $A$ individual of age $s$ in the homogeneous process has $\text{Poisson}( s \, r^*_A)$ offspring we find
\[
\pr(v \text{ is of type A and has }k \text{ offspring})
\, \approx \,
p \int_0^{\infty} \pr(\text{Poisson}( s \, r^*_A)=k ) \, e^{-s} \, \mathrm{d}s
\, = \, 
\frac{p}{1+r^*_A} \left( \frac{r^*_A}{1+r^*_A} \right)^k.
\] 
A similar heuristic calculation can be used to motivate the limiting degree distribution for type $B$ individuals.
Note here however, that in our formal proofs we will not apply results from Jagers and Nerman \cite{jagers1996asymptotic} since we have {\bf pdBP} instead of time homogeneous multi-type branching process.

It is worth noting that one may also carry out an approximation like that above for the URT.  
In this case the corresponding process is a (unitype) continuous time branching 
process with Poisson offspring rate 1. 
Results from Jagers and Nerman \cite{jagers1996asymptotic} show that, for the branching process,
the limiting age distribution for a randomly chosen individual $v$ is Exponential(1) and we find
\[
\pr(v \text{ has } k \text{ offspring}) \, \approx \,
\int_0^{\infty} \pr(\text{Poisson}(s)=k) \, e^{-s} \, \mathrm{d}s 
\, = \,
2^{-k-1},
\] 
which is exactly the limiting degree distribution of URT.
   
Now suppose $p=1$, using the construction in special case \ref{item:spec-construct} of Section \ref{sec:special} we have that for a randomly chosen individual $v$ of type $A$ in $\BP(t)$  
\[
\pr(v \text{ has } k \text{ offsprings}) \approx
\sum_{i=0}^{\infty} 2^{-k-i-1} \binom{k+i}{k} q^k (1-q)^i = 
\frac{ q^k }{ 2^{k+1} } \sum_{i=0}^{\infty} \binom{k+i}{k} \left( \frac{1-q}{2} \right)^i = 
\frac{1}{q+1} \left( \frac{q}{q+1} \right)^k.
\]
Here the infinite sum is calculated using binomial series, which is equivalent to summing up the probability mass function of a negative binomial distribution in this case.

\vskip.1in

Using the limiting degree distribution, one may derive consistent estimators of the model parameters $p$ and $q$ using
the observed degrees of the tree.  It is sufficient to consider the statistics $N_0(n)$, the number of leaves of $\cT_n$, 
and $N_1(n)$, the number of vertices with out-degree 1 in $\cT_n$.

\begin{cor}
	\label{cor:estimators}
If $p \neq 1/2$ and $q \neq 1$ then there exist consistent estimators $\hat{p}_n$ and $\hat{q}_n$ for $p$ and $q$ 
that can be computed by solving a quadratic equation involving only $N_0(n)$ and $N_1(n)$ (as given by \eqref{eqn:estimator-1} and \eqref{eqn:estimator-2}).
\end{cor}

\begin{rem}
In practice, if the number of vertices $n$ is large, one may use a sub-sample to estimate $N_0(n)$ and $N_1(n)$. 
If these estimates are 
consistent, the estimates of $p$ and $q$ obtained from them will be consistent as well.
\end{rem}

\begin{rem}
As noted above, when $p=1/2$ or $q=1$ the limiting degree distribution if the CMRT matches that of the URT. 
In this setting it follows from results of \cite{aldous1991asymptotic} and \cite{jagers1996asymptotic} 
that the \emph{local weak limits} of both the URT and the CMRT exist and are equal to the same infinite {\tt sin}-tree. As such, distinguishing these two models using the densities of local statistics appears to be a difficult problem. 
Still, global statistics might remain informative. For example, when $q=1$, the CMRT is just two disjoint URTs with the two roots connected. Using similar algorithm as introduced by Bubeck, Devroye and Lugosi \cite{bubeck2017finding} one can construct a confidence set for the roots and estimate $p$ from that.
\end{rem}

The next two results concern the limiting behavior of two global statistics for the CMRT $\cT_n$, 
specifically, the maximal degree $M_n$ of any vertex, and 
the height $H_n$ of the tree, which we define to be the maximum distance from any vertex 
to the nearest root. 

\vskip.1in

\begin{thm}
	\label{thm:max-degree}
	Let $M_n$ denote the maximal degree of any vertex in the CMR $\cT_n$. 
	When $p,q \neq 1$, there exist constants $C_1, C_2 > 0$ depending only on $p$ and $q$ such that
	\[ 
	\liminf_{n \to \infty}  \frac{ M_n }{\log n} \geq C_1
	\quad \text{and} \quad
	\limsup_{n \to \infty}  \frac{ M_n }{\log n} \leq C_2 
	\quad a.s. 
	\]
	Moreover, when $q=1$, 
	\[
	\frac{ M_n }{\log n} \convas \frac{1}{\log 2},
	\]
	and when $p=1$ but $q \neq 1$,
	\[
	\frac{ M_n }{n} \convas 1-q.
	\]
\end{thm}

\vskip.3in

\begin{thm}
	\label{thm:height}
	Let $H_n$ denote the height of the CMR $\cT_n$. When $p \neq 1$ or $q=1$,
	\[
	\frac{H_n}{\log n} \convas e.
	\]
	When $p=1$ and $q \neq 1$,
	\[ 
	\liminf_{n \to \infty}  \frac{ H_n }{\log n} \geq qe 
	\quad \text{and} \quad
	\limsup_{n \to \infty}  \frac{ H_n }{\log n} \leq e 
	\quad a.s. 
	\]
\end{thm}

\vskip.2in

\begin{rem}
When $p=1$ and $q=0$, the tree $\cT_n$ looks like a ``star''  in which $n-1$ type $A$ vertices are 
connect to the type $B$ root. In this case the height of the tree is $H_n = 1$, and 
$\lim\limits_{n \to \infty}  H_n / \log n = 0$.
\end{rem}

\vskip.2in
	
\subsection{Results for general CMRT}
\label{sec:cmrt-res-general}
We now present results for $K$ type CMRTs with arbitrary attachment probabilities.

\vskip.1in

\begin{thm}
	\label{thm:degree-dist-general}
	For each fixed $k$, let $N_k(n)$ denote the number of vertices with out-degree $k$ in a 
	$K$ type CMRT $\cT_n$. Then 
	\[
	\frac{N_k(n)}{n} \convp c_k
	\] 
	where 
	\begin{equation*}
		c_k := \sum_{i=1}^K \frac{p_i}{ 1 + r_i } \left( \frac{r_i}{ 1 + r_i } \right)^k.
	\end{equation*}
	Here 
	\[
	r_i = \frac{1}{p_i} \sum_{j=1}^K p_j q_{ji}. 
	\]
	When $p_i=0$ the $i$-th term in $c_k$ should be interpreted as 0.
\end{thm}

\vskip.1in

\begin{rem}
Analogously to the CMRT with two types, the limiting degree distribution of a $K$-type CMRT is a mixture of $K$ shifted 
geometric distributions, which coincides with that of a URT if and only if the balance equation 
$\sum_{j=1}^K p_j q_{ji} = p_i$ holds for all $1 \leq i \leq K$. Note that $\sum_{j=1}^K p_j q_{ji}$ 
is the probability that the parent of a new vertex is of type $i$. 
Thus the balance equation essentially states that the type distribution for parents is the same as that of 
their children.  In the continuous time embedding, the condition implies that all individuals will reproduce 
at approximately the same rate once the population stabilizes. 
\end{rem}

\vskip.1in

\begin{thm}
	\label{thm:max-degree-general}
	Let $M_n$ denote the maximal degree in the $K$-type CMRT $\cT_n$. 
	When $p_1 < 1$ and $q_{ii} < 1$ for some $1 \leq i \leq K$ there exists constants $C_1$ and $C_2$ 
	depending only on $p_i$ and $q_{ij}$ such that
	\[ 
	\liminf_{n \to \infty}  \frac{ M_n }{\log n} \geq C_1
	\quad \text{and} \quad
	\limsup_{n \to \infty}  \frac{ M_n }{\log n} \leq C_2
	\quad a.s. 
	\]
	When $q_{ii} = 1$ for all $1 \leq i \leq K$, 
	\[
	\frac{ M_n }{\log n} \convas \frac{1}{\log 2}
	\]
	and when $p_1=1$ but $q_{11} < 1$,
	\[
	\frac{ M_n }{n} \convas \max_{2 \leq i \leq K} q_{1i}.
	\]
\end{thm}

\vskip.3in

\begin{thm}
	\label{thm:height-general}
	Let $H_n$ denote the height of the $K$-type CMRT $\cT_n$.  When $p_1 < 1$ or $q_{11}=1$,
	\[
	\frac{H_n}{\log n} \convas e.
	\]
	When $p_1=1$ and $q_{11} < 1$,
	\[ 
	\liminf_{n \to \infty}  \frac{ H_n }{\log n} \geq q_{11} e 
	\quad \text{and} \quad
	\limsup_{n \to \infty}  \frac{ H_n }{\log n} \leq e 
	\quad a.s. 
	\]
\end{thm}

\vskip.2in

\subsection{Limiting degree distribution of CWRT}
Following the same arguments as in the case of community modulated recursive trees we can derive the limiting degree 
distribution of $K$-type community weighted recursive trees.

\vskip.1in

\begin{thm}
	\label{thm:degree-dist-cwrt}
	For each fixed $k$, let $N_k(n)$ denote the number of vertices with out-degree $k$ in
	a $K$-type CWRT $\tilde{\cT}_n$. Then 
	\[
	\frac{N_k(n)}{n} \convp \tilde{c}_k
	\] 
	where 
	\begin{equation*}
		\tilde{c}_k := \sum_{i=1}^K \frac{p_i}{ 1 + \tilde{r}_i } \left( \frac{\tilde{r}_i}{ 1 +  \tilde{r}_i } \right)^k
	\end{equation*}
	with
	\[
	\tilde{r}_i = \sum_{j=1}^K \frac{p_j \omega_{ji}}{\sum_{l=1}^K p_l \omega_{jl}}. 
	\]
\end{thm}

\vskip.1in

\begin{rem}
Analogously to CMRTs, the limiting degree distribution of a $K$-type CWRT is a mixture of $K$ geometric distributions
with different parameters. 
In the special case $\omega_{ij} \equiv 1$ the limiting distribution is identical to that of URT.
\end{rem}

\vskip.2in

\subsection{Limiting degree distribution of CMPA}
We now describe results for the community modulated preferential attachment. Let $p^*_i = \sum_{j=1}^K p_j q_{ji}$ be the probability for the parent of a new vertex to be of type $i$ in CMPA (same as in CMRT). We also derive the limiting degree distribution of CMPA:

\vskip.1in

\begin{thm}
	\label{thm:degree-dist-cmpa}
	For each fixed $k$, let $N_k(n)$ denote the number of vertices with out-degree $k$ in $\cT^*_n$. Then 
	\[
	\frac{N_k(n)}{n} \convp c^*_k
	\] 
	where 
	\begin{equation}
	\label{const-cmpa}
	c^*_k := \sum_{i=1}^K \frac{p_i}{\nu_i}  \frac{\Gamma(k+\alpha_i) \Gamma(1/\nu_i + \alpha_i)}{\Gamma(\alpha_i) \Gamma(k+1+1/\nu_i + \alpha_i)}.
	\end{equation}
	Here 
	\[
	\nu_i = \sum_{j=1}^K \nu_{ji}
	\]
	with
	\[
	\nu_{ji} = \frac{p_j q_{ji}}{\alpha_{ji} p_i + p^*_i},
	\]
	and
	\[
	\alpha_i =  \frac{\sum_{j=1}^K \nu_{ji} \alpha_{ji}}{\sum_{j=1}^K \nu_{ji}} = 
	\frac{\sum_{j=1}^K \nu_{ji} \alpha_{ji}}{\nu_i}.
	\]
	And $\Gamma(\cdot)$ is the gamma function. When $p^*_i =0$, the $i$-th term in $c^*_k$ should be interpreted as $p_i$ when $k=0$, and 0 otherwise.
\end{thm}

\vskip.1in

\begin{rem}
	Similar to CMRT, the limiting degree distribution of CMPA tree is a mixture of $K$ distributions, and coincides with that of preferential attachment model \cite{bollobas2001degree} when $\alpha_{ji} \equiv \alpha$ and the balance equation $p^*_i = \sum_{j=1}^K p_j q_{ji} = p_i$ holds for all $1 \leq i \leq K$. Moreover, for large $k$, using Stirling's approximation we get $c^*_k = O(k^{-1-1/\nu})$ where $\nu = \max_i \nu_i$. Consider the special case where $\alpha_{ji} \equiv \alpha$. we have $\sum_{i=1}^K p^*_i = \sum_{j=1}^K p_j \sum_{i=1}^K q_{ji} = \sum_{j=1}^K p_j = 1$, and 
	$\sum_{i=1}^K p^*_i (1/\nu_i -1) = \alpha \sum_{i=1}^K p_i = \alpha$. 
	It follows that $-1 - 1/\nu \geq -2-\alpha$ and equality holds if and only if the balance equation holds. Therefore, the limiting degree distribution of CMPA tree typically has a heavier tail than that of preferential attachment tree with the same parameter $\alpha$. The tail is heavier when $p^*_i$ is larger compared to $p_i$, which agrees with intuition. 
\end{rem}

\vskip.2in

\section{Related work and discussion} 
\label{sec:related}

\begin{enumeratea}
\item {\bf Branching processes and growing trees:} Branching processes have been used to study a variety of tree based stochastic structures. See Devroye \cite{devroye1998branching} 
for a nice summary of classical results that can be proved via
branching processes, both discrete and continuous. 
In this paper we make use of continuous time branching processes. 
Since Athreya and Karlin \cite{athreya1968embedding} introduced the continuous time embedding of urn processes, 
embeddings of this sort has been used to study a variety of tree structures. 
The well-studied Crump-Mode-Jagers branching
process \cite{crump1968general} is known to give rise to a rich family of stochastic trees, of which uniform recursive trees are a special case.  Indeed, some of the results for the URT described below were proved via continuous time embeddings \cite{devroye1995strong}\cite{pittel1994note}. 
Although the standard techniques for analyzing asymptotic degree distributions were developed under the discrete setting, it could be fruitful to work with continuous time embeddings for more complicated models \cite{bhamidi2015change}. 
In general, this technique is more mathematically involved, but when tractable allows access to derive asymptotics for a host of functionals including degree distribution, maximal degree height and so on. 
	
\item {\bf Related work on the URT:} Our results on the asymptotics of CMRT are closely related to that of URT, which have been studied extensively. 
Let $\set{\cU_n}_{n\geq 1}$ be a URT, and $N_k(n)\ (k \geq 0)$ denote the number of vertices with out-degree $k$ in $\cU_n$. 
Moon\cite{moon1974distance} showed that 
\[
\frac{N_k(n)}{n} \convp 2^{-k-1},
\]  
and Janson\cite{janson2005asymptotic} further showed that 
\[
n^{-\frac{1}{2}}(N_k(n) - 2^{-k-1} n) \convd V_k
\] 
jointly for all $k\geq 0$, where  $\set{V_k}_{k\geq 0}$ is a Gaussian process.

Let $M_n$ and $H_n$ denote respectively the maximal degree and height in $\cU_n$. Here height is defined as length of the longest upward path from a leaf to the root. Devroye and Lu\cite{devroye1995strong} showed that 
\[
\frac{M_n}{\log_2 n} \convas 1 
\quad \text{and} \quad 
\lim\limits_{n \to \infty}\frac{\rE M_n}{\log_2 n}=1.
\] 
And Pittel\cite{pittel1994note} showed that 
\[
\frac{H_n}{\log n} \convas e.
\]

\item {\bf Related work on the preferential attachment tree:} The preferential attachment tree and its variants have also been extensively studied in the literature. Recall that in linear preferential attachment trees, instead of uniform attachment, a new vertex connects to an existing vertex with probability proportional to $\alpha$ plus that vertex's out-degree. The linear preferential attachment tree is known to have a limiting 
power law distribution with exponent depending on the model parameter $\alpha$ \cite{bollobas2001degree}. 
In terms of models with multiple types, which are our focus here, to the best of our knowledge the only result we are aware of is Deijfen and Fitzner \cite{deijfen2016birds}, who computed the limiting degree distribution for a special case of CMPA with two types heuristically and conducted a simulation study.

\item {\bf Community detection:}  
As described in the introduction, one motivation for this work was 
to draw connections between the research area of community detection and that of stochastic trees. 
It is well nigh impossible to give a representative set of references to community detection, but the \cite{karrer2011stochastic,newman2012communities,fortunato2010community} and the references 
therein provide a good overview of applications, while \cite{bordenave2015non,krzakala2013spectral} 
provide an introduction to the burgeoning literature in the probability community. 
Loosely speaking, community detection is an unsupervised learning task (closely related to clustering) for obtaining insight into networks, in particular understanding subsets within a network, vertices within which seem more densely connected within the subset as opposed to outside the subset. A host of techniques have been proposed to extract such subsets within networks 
and there has also been an enormous effort in the probability community to evaluate the performance of proposed techniques on tractable network models with known community structure. The aim of this work was to derive math results for the asymptotics of various functionals related to community modulated growth schemes in networks. In work in progress we propose and study the performance of various techniques for extracting the underlying latent community structure if we only have access to the graph without any information on the corresponding types of vertices. 
 
\end{enumeratea}

\section{Proofs}
\label{sec:proof}
This section is organized as follows. First in Section \ref{sec:proof-em} we prove the equivalence of the continuous time embedding, together with some of its basic properties. Then in Section \ref{sec:proof-deg} we derive the limiting degree distribution for Community Modulated Recursive Tree (CMRT) with two types and use it in Section \ref{sec:proof-est} to give consistent estimators of model parameters. In Section \ref{sec:proof-max-deg} and \ref{sec:proof-height} we prove results for maximal degree and height. Finally, Section \ref{sec:proof-cmrt-general}, \ref{sec:proof-cwrt} and \ref{sec:proof-cmpa} extend the results to general CMRT and the other two variants: Community Weighted Recursive Tree (CWRT) and Community Modulated Preferential Attachment (CMPA).

\subsection{Proof for continuous-time embedding and some basic properties}
\label{sec:proof-em}
In this section we shall first give a proof of the continuous-time embedding (Lemma \ref{lem:embed}).

{\bf Proof:} Assume that 
$\set{\BP(T_n)}_{2 \leq n \leq k}  \stackrel{d}{=} \set{\cT_n}_{2 \leq n \leq k}$ 
for a fixed integer $k\geq 2$ (which holds for $k=2$ by definition). Conditioning on $\set{\BP(T_n)}_{2\leq n\leq k}$, it can be checked using properties of exponential distribution that the probability for the next individual born to be of type $A$ is 
\[
\frac{ n_A(T_n) r_{AA}(T_n) + n_B(T_n) r_{BA}(T_n) }
{ n_A(T_n) (r_{AA}(T_n) + r_{AB}(T_n)) + 
	n_B(T_n) (r_{BA} (T_n) + r_{BB} (T_n)) } = p
\] 
and the probability for the next type $A$ individual born to have a parent of type $A$ is
\[
\frac { n_A(T_n) r_{AA}(T_n) } 
{ n_A(T_n) r_{AA}(T_n) + n_B(T_n) r_{BA}(T_n) } = q.
\] 
Also, given a type $\in \set{A,B}$, the probability for each individual of this type to give birth to the next type $A$ individual is equal. Similarly one can check the corresponding probability for type $B$ individuals. Thus the dynamics of $\BP(T_{k+1})$ conditioning on $\set{\BP(T_n)}_{2\leq n\leq k}$ is the same as that of $\cT_{k+1}$ conditioning on $\set{\cT_n}_{2\leq n\leq k}$. Therefore by induction we have the desired result.

\qed

Now that we have the continuous-time embedding, we shall proceed to derive some of its basic properties that will come in handy. In what follows, we shall assume that natural filtration $\set{\cF(t) : t \geq 0}$ is used throughout.

\begin{lem}\label{lem:l2-mart}
	The process $\set{e^{-t} n_A(t)}_{t\geq 0}$ is an $\bL^2$-bounded positive martingale. In particular there exists a strictly positive finite random variable $W$ such that 
	\[
	e^{-t} n_A(t) \convas W, \qquad \mbox{ as } t \to \infty. 
	\]
\end{lem}
\begin{rem}
	As will be evident from the calculation below, the marginal distribution of $n_A(\cdot)$ is identical to that of a rate one Yule process starting with a single individual. In particular the limit random variable $W \stackrel{d}{=} \exp(1)$. Recall that a Yule process with rate $\lambda$ is a time-inhomogeneous Poisson process with birth rate $\lambda i$, where $i$ is the current population size.
\end{rem}

{\bf Proof:} First we introduce some preliminary notations that will be used extensively throughout the proofs. Recall that for a jump diffusion 
$\set{ X(t) \in \bR^n }_{ t\geq 0 }$, 
its infinitesimal generator $\cA$ is defined for functions $f: \bR^n\to\bR$ by 
\[
\cA f(x) = \lim\limits_{t\to 0^+} \frac{1}{t} ( \bE( f(X(t)) | X(0)=x ) - f(x) )
\] 
if the limit exists. Then by Dynkin's formula (see {\O}ksendal and Sulem \cite[Chapter 1.3]{oksendal2003stochastic} for a formulation) and Markov property of jump diffusion we have that 
$\set{ X(t) - \int_0^t \cA X(s)ds }_{ t \geq 0 }$ 
is a martingale. Note that if the diffusion term is not present and the jump part is an inhomogeneous Poisson process, we have that 
$\cA X(t) = \delta(X(t)) \lambda(X(t))$ 
where $\delta(x)$ and $\lambda(x)$ are size and intensity of jump when the process is at $x\in\bR^n$.

Now we are ready to introduce some martingales. Denote $\exp(-t) n_A(t)$ by $\tilde{n}_A(t)$. 
Using the rates in Section \ref{sec:ct-em} we have 
\[
\cA \tilde{n}_A(t) = 
e^{-t} n_A(t) - e^{-t} n_A(t)=0
\] 
and 
\[
\cA\tilde{n}_A^2(t) = 
e^{-2t} (2n_A(t) + 1) n_A(t) - 2 e^{-2t} n_A^2(t) = 
e^{-2t} n_A(t).\]

Thus by Dynkin's formula we have that both $\set{ \tilde{n}_A(t) : t \geq 0 }$ and 
\[
\tilde{n}^2_A(t) - \int_0^t e^{-2s} n_A(s) ds, 
\qquad 
t\geq 0
\] 
are martingales. Taking expectation of both martingales we get 
$\E(n_A(t)) = \exp(t)$ and $\E(\tilde{n}^2_A(t)) = 2-\exp(-t)$. Therefore $\set{\tilde{n}_A(t):t\geq 0}$ is $\bL^2$ bounded and the second statement follows from standard martingale convergence theorem.

\qed

\begin{lem}
	\label{lem:difference}
	Define $Z(t) := pn_B(t) - (1-p)n_A(t)$. Then $\set{Z(t):t\geq 0}$ is a martingale and further $e^{-t} Z(t) \convas 0$. This implies
	$e^{-t} n_B(t) \convas (1-p) W/p$ 
	where $W$ is as in Lemma \ref{lem:l2-mart}.  
\end{lem}
{\bf Proof:} 
Again using the rates in Section \ref{sec:ct-em} we have 
\[
\cA Z(t) = 
p \frac{1-p}{p} n_A(t) - (1-p) n_A(t) = 0
\] 
and 
\[
\cA Z^2(t) =
p^2 \cA n_B^2(t) + (1-p)^2 \cA n_A^2(t) - 2p(1-p) \cA n_A n_B(t),
\] 
with 
\[
\cA n_B^2(t) = 
(2n_B(t)+1) \frac{1-p}{p} n_A(t), 
\quad 
\cA n_A^2(t) = 
(2n_A(t)+1) n_A(t)
\]
and 
\[
\cA n_A n_B(t)= 
n_A(t) \frac{1-p}{p} n_A(t) + n_B(t) n_A(t).
\]

Thus by some elementary algebra and Dynkin's formula we have that both $\set{Z(t):t\geq 0}$ and
\[
M(t) = 
Z^2(t) - \int_0^t (1-p) n_A(s) ds, 
\qquad t \geq 0
\]
are martingales. From Lemma \ref{lem:l2-mart} we get that 
$\E(n_A(s)) = \exp(s)$. Taking expectation of $\E(M(t))$ shows 
\[
\E(Z^2(t)) = 
(1-p) (e^t -1) + (2p-1)^2. 
\]

Now apply Markov's inequality to $\exp(-4\log n) Z^2(2\log n)$ we have for any $\epsilon>0$:
\[
\pr \left( e^{-4\log n} Z^2(2\log n) > \epsilon \right) \leq 
\frac{ e^{-4\log n }( (1-p) (e^{2\log n} -1) + (2p-1)^2 ) }{\epsilon} = 
\frac{1-p}{ n^2 \epsilon } + \frac{ 4p^2 - 3p }{ n^4\epsilon }.
\]
Thus by the Borel--Cantelli lemma $\exp(-2\log n) Z(2\log n) \convas 0$. 
Since we know by Lemma \ref{lem:l2-mart} that 
$\exp(-2\log n) n_A(2\log n) \convas W$, 
we get $\exp(-2\log n) n_B(2\log n) \convas (1-p) W/p$.

Finally for any $t>0$ we can find a positive integer $n$ such that 
$2\log n \leq t < 2\log(n+1)$. By monotonicity of $n_B(\cdot)$ we have
$n_B(2\log n) \leq n_B(t) < n_B(2\log(n+1))$ and further 
\[
e^{-2\log(n+1)} n_B(2\log n) \leq 
e^{-t} n_B(t) < 
e^{-2\log n} n_B(2\log(n+1)).
\] 
Since the left hand side 
\[
e^{-2\log(n+1)} n_B(2\log n) = 
e^{-2\log n} n_B(2\log n) \cdot \frac{n^2}{(n+1)^2} 
\convas 
\frac{ (1-p)W }{p}
\] 
and similarly the right hand side converges a.s. to the same limit, we have $\exp(-t) n_B(t) \convas (1-p) W/p$. This immediately implies $\exp(-t) Z(t) \convas 0$.

\qed
\begin{rem}
	This lemma essentially proves that 
	$n_B(t)/n_A(t) \convas (1-p)/p$, which is used in Section \ref{sec:cmrt-res} for heuristic calculations. 
\end{rem}

\begin{lem}
	\label{lem:stop-time}
	The population size process $n(t)$ satisfies 
	$e^{-t} n(t)\convas W + (1-p) W/p := W_\infty $. 
	In particular, the sequence of stopping times $T_n$ satisfy
	\[
	T_n - \log{n} \convas -\log(W_\infty),
	\]
	where $W_\infty$ is a strictly positive finite random variable. 
\end{lem}
The first statement here is a direct corollary of Lemma \ref{lem:l2-mart} and \ref{lem:difference}, and the second statement follows by replacing $t$ with $T_n$. Note that the first statement essentially says that the Malthusian parameter for process $\set{ n(t) : t \geq 0 }$ is 1.

\subsection{Proof for limiting degree distribution}
\label{sec:proof-deg}
In this section we shall prove Theorem \ref{thm:degree-dist}.

To work in continuous time, we need to reformulate Theorem \ref{thm:degree-dist} via the embedding. Let $N_{k,A}(n)$ and $N_{k,B}(n)$ denote the number of type $A$ and $B$ vertices with out-degree $k$ in $\cT_n$.  The plan is as follows: first we focus on type $A$ vertices and prove the result below, and as similar result holds for type $B$ vertices, using $N_k(n) = N_{k,A}(n) + N_{k,B}(n)$ we completes the proof. Note that when $p=1$, there is no need to consider type $B$ vertices.

\begin{thm}
	\label{thm: A-deg-dist}
	For a fixed integer $k>0$, let $n_{k,A}(t)$ denote the number of type $A$ individuals with $k$ offsprings in $\BP(t)$. Then
	\[
	\frac{n_{k,A}(t)} {n(t)} 
	\convp p_{k,A}, 
	\qquad 
	\text{as } t \to \infty. 
	\]
	Here
	\begin{equation*}
	\label{eqn:pka-def}
	p_{k,A} := 
	p \int_0^{\infty} \pr(\mathrm{Poisson} (r^*_A s)=k) e^{-s} ds = 
	\frac{p}{1+r^*_A} \left( \frac{r^*_A}{1+r^*_A} \right)^k
	\end{equation*}
	where $r^*_A$ is the total reproduction rate of type A individuals:
	\begin{equation}
	\label{eqn:ra-def}
	r^*_A = q + \frac{ (1-p)(1-q) }{p}. 
	\end{equation} 
	
	Thus by the embedding in Lemma \ref{lem:embed} for $\set{\cT_n}_{n\geq 2}$ we have
	\[
	\frac{ N_{k,A}(n) }{n} \convp p_{k,A}, 
	\qquad 
	\text{as } n \to \infty. 
	\]
\end{thm}

{\bf Proof:} 
Throughout the proof we work with the continuous time embedding. For any fixed constant $0<a<t$, let $n_{k,A}[t-a,t]$ be the number of type $A$ individuals born in the interval $[t-a,t]$ that have exactly $k$ offsprings by time $t$. Given Lemma \ref{lem:stop-time}, it is enough to show the following two propositions.
\begin{prop}
	\label{prop:trunc-diff}
	\[
	\limsup_{a \to \infty} \limsup_{t \to \infty} 
	e^{-t} \left( n_{k,A}(t) - n_{k,A}[t-a,t] \right) = 0,
	\qquad a.s. 
	\]
\end{prop}
{\bf Proof:} Since the population size $n_A(t)$ grows exponentially, most type $A$ individuals are born after time $t-a$. Indeed, $n_{k,A}(t) - n_{k,A}[t-a,t] = n_{k,A}(t-a) \leq n_A(t-a)$, and by Lemma \ref{lem:l2-mart} we have 
\[
\limsup_{t \to \infty} 
e^{-t} \left( n_{k,A}(t) - n_{k,A}[t-a,t] \right) \leq
e^{-a} \lim_{t \to \infty} 
e^{-(t-a)}n_A(t-a) = 
e^{-a} W,
\qquad a.s. 
\] 
Letting $a \to \infty$ proves the proposition.

\qed

\begin{prop}
	\label{prop:trunc-deg-dist}
	Recall the random variable $W$ in Lemma \ref{lem:difference}. 
	Then for each fixed $a>0$,  we have
	\[
	e^{-t} n_{k,A}[t-a,t] \convp 
	W \int_0^a \pr(\mathrm{Poisson}(r^*_A s)=k) e^{-s} ds 
	\]
	as $t \to \infty$.
\end{prop}

This assertion needs some work and the proof follows a similar procedure as that in Bhamidi et al. \cite[Section 4.2.2]{bhamidi2015change}. First, recall from Lemma \ref{lem:l2-mart} that $n_A(t) \approx W e^t$ for large $t$. For our proof we will need a finer concentration result that goes as follows:

\begin{lem}\label{lem:na-sup}
	\[
	\pr \left( \sup_{t-a \leq s \leq t} |n_A(s) - W e^s| < \sqrt{t e^t} \right) \to 1
	\]
	as $t \to \infty$ where $W$ is as in Lemma \ref{lem:l2-mart}. Equivalently, we shall say that w.h.p. as $t \to \infty$,
	$\sup_{t-a \leq s \leq t} |n_A(s) - W \exp(s)| < \sqrt{t \exp(t)}$.
\end{lem}
{\bf Proof:} First note that 
\[
e^{-s} n_A(s) - e^{-(t-a)} n_A(t-a),
\qquad s \geq t-a
\]
is a martingale by Lemma \ref{lem:l2-mart} and recall that 
$\E (\exp (-2s) n_A^2(s) ) = 2 - \exp(-s)$ from the proof.

Now fix any $T > t-a$, Doob's $L^2$-maximal inequality (here we use the version stated in {\O}ksendal \cite[Theorem 3.2.4]{oksendal2003stochastic}, see also Karatzas and Shreve \cite[Theorem 3.8]{karatzas2012brownian} for a proof) applied to the above martingale gives, for any $C>0$:

\begin{equation}
\label{eqn:na2-doob-1}
\pr \left( \sup_{t-a \leq s \leq T} |e^{-s} n_A(s) - e^{-(t-a)} n_A(t-a)| > C \right)
\leq \frac{\E (e^{-T} n_A(T) - e^{-(t-a)} n_A(t-a))^2}{C^2}
\end{equation}

Recall that $\set{\tilde{n}_A(t) := \exp(-t) n_A(t)}_{t\geq 0}$ is a martingale so we have
\[
\E (\tilde{n}_A(T) \tilde{n}_A(t-a)) =
\E \left( \tilde{n}_A(t-a)
\E \left( \tilde{n}_A(T) | \tilde{n}_A(t-a) \right)  
\right) =
\E ( \tilde{n}_A^2(t-a) ) =
\E (e^{-2(t-a)} n_A^2(t-a)).
\]
It follows that
\[
\E (e^{-T} n_A(T) - e^{-(t-a)} n_A(t-a))^2 = 
\E (e^{-2(t-a)} n_A^2(t-a)) - \E (e^{-2T} n_A^2(T)) =
e^{-(t-a)}-e^{-T}.
\]
Plug this in \eqref{eqn:na2-doob-1} we get 
\begin{equation}
\label{eqn:na2-doob-2}
\pr \left( \sup_{t-a \leq s \leq T} |e^{-s} n_A(s) - e^{-(t-a)} n_A(t-a)| > C \right)
\leq \frac{ e^{-(t-a)}-e^{-T} }{C^2}
\end{equation}

Now, let $T \to \infty$ and use the a.s. convergence result from Lemma \ref{lem:l2-mart} to yield 

\begin{equation*}
\label{eqn:na2-doob-3}
\pr \left( |W - e^{-(t-a)} n_A(t-a)| > C \right) \leq 
\frac{ e^{-(t-a)} }{C^2}.
\end{equation*}

On the other hand, let $T=t$ in \eqref{eqn:na2-doob-2} we get

\begin{equation*}
\label{eqn:na2-doob-4}
\pr \left( \sup_{t-a \leq s \leq t} |e^{-s} n_A(s) - e^{-(t-a)} n_A(t-a)| > C \right)
\leq \frac{ e^{-(t-a)} - e^{-t} }{C^2} 
\leq \frac{ e^{-(t-a)} }{C^2}.
\end{equation*}

Finally, combine the above two inequalities we have
\begin{equation*}
\label{eqn:na2-doob-5}
\pr \left( \sup_{t-a \leq s \leq t} |e^{-t} n_A(t) - W| > 2C \right) \leq
\frac{ 2e^{-(t-a)} }{C^2}.
\end{equation*}

Let $C = \sqrt{t\exp(-t)} / 2$ and after some simple algebraic manipulation we get the wanted result.

\qed

Now that we have Lemma \ref{lem:na-sup}, we will proceed to approximate the integral in Proposition \ref{prop:trunc-deg-dist} by $e^{-t} n_{k,A}[t-a,t]$. To do so, we divide the interval $[t-a,t]$ into intervals of length $\delta:=e^{-t/3}$ and denote by $\cI_i (1\leq i\leq a/\delta)$ the $i$-th interval $[t-a +(i-1)\delta, t-a+i\delta]$ (for simplicity we treat $a / \delta$ as if it is an integer but this should not matter). Let $t_i=t-a +(i-1)\delta$ and write $n_A(\cI_i)$ for the number of type $A$ individuals born in $\cI_i$.

On the event given in Lemma \ref{lem:na-sup}, we have
\[
|n_A( t_i + \delta ) - n_A(t_i) - W e^{t_i} (e^\delta-1) | < 
2 \sqrt{t e^t}.
\]
Since
$W \exp(t_i) (\exp(\delta)-1) = 
\delta W \exp(t_i) + o_p ( \sqrt{\exp(t)} )$, we can further get
\begin{equation}
\label{eqn:na-uniform}
\pr \left( \bigcap^{a/\delta}_{i=1} \Big\{ |n_A(\cI_i) - \delta W e^{t_i}| < 
3 \sqrt{t e^t} \Big\} \right) \to 1
\end{equation}
as $t \to \infty$. This shows that the number of type $A$ individuals born in $\cI_i$ is approximately $\delta W e^{t_i}$.

Unfortunately, even individuals born in the same interval could start reproducing at different times (\ie have distinct birth times), and we need to at least align these individuals. To facilitate the analysis, for each type $A$ individual born in $\cI_i$, we call it ``good'' if it gives birth to no offspring (type $A$ or $B$) in $\cI_i$ and ``bad'' otherwise. Since the interval is small, we expect most individuals born to be good:

\begin{lem}\label{lem:bad}
	There exists a constant $M>0$ such that 
	\[
	\pr \left( \bigcap^{a / \delta}_{i=1} \Big \{ n_A^{bad}(\cI_i) < 
	M W t e^{t/3} \Big \} \right) 
	\to 1
	\]
	as $t \to \infty$, where $n_A^{bad}(\cI_i)$ is the number of bad individuals born in $\cI_i$.
\end{lem}

{\bf Proof:} Call a bad individual in $\cI_i$ a ``direct'' bad individual if it is an offspring of individuals born before $\cI_i$ and write $n_A^{dir}(\cI_i)$ for the number of direct bad individuals born in $\cI_i$. As non-direct bad individuals in $\cI_i$ must have an antecedent that is direct bad in $\cI_i$, we have $n_A^{bad}(\cI_i) \leq n_A^{desc}(\cI_i)$, where $n_A^{desc}(\cI_i)$ is the number of descendants of direct bad individuals in $\cI_i$.

Note that a direct bad individual in $\cI_i$ has to satisfy the following two conditions:
\begin{itemize}
	\item First it has to be an offspring of individuals born before $\cI_i$ and the number of such offsprings is at most $n_A(\cI_i)$.  
	\item Secondly it has to give birth to at least one offspring in $\cI_i$ and the corresponding (conditional) probability is at most 
	\begin{equation*}
	\label{eqn:pbad-def}
	p_{bad} := \pr (\mathrm{Exp} (r^*_A) \leq \delta) \sim r^*_A \delta.
	\end{equation*}
\end{itemize}
Combining the above conditions we have 

\begin{equation*}
\label{eqn:ndir-bound}
n_A^{dir} (\cI_i) \preceq \mathrm{Bin} (n_A(\cI_i) , p_{bad} ).
\end{equation*}
where $\preceq$ denotes stochastic dominance. By \eqref{eqn:na-uniform} we have that w.h.p. as $t \to \infty$, 
$W \exp(2t/3) / 2 \leq n_A(\cI_i) \leq 2W \exp(2t/3)$ for all $i$. Condition on this event and use Hoeffding inequality for binomial distribution we get
\[
\pr \left( n_A^{dir}(\cI_i) > 
n_{\delta} p_{bad} + \sqrt{n_{\delta}\log(n_{\delta})} \right) \leq
\frac{1}{n_A^2(\cI_i)} \leq
\frac{4}{ W^2 \exp(4t/3) }
\]
where $n_{\delta} := 2W \exp(2t/3)$.
Since there exists $C_1>0$ such that 
$n_{\delta} p_{bad} + \sqrt{n_{\delta}\log(n_{\delta})} < C_1 W t e^{t/3}$
for large enough $t$, and 
\[
\frac{a}{\delta} \cdot \frac{4}{ W^2 \exp(4t/3) } \to 0
\]
as $t \to \infty$, using a union bound we have
\begin{equation}
\label{eqn:ndir-uniform}
\pr \left( \bigcap^{ a / \delta }_{i=1} 
\Big\{ n_A^{dir}(\cI_i) \leq C_1 W t e^{t/3} \Big\} \right) 
\to 1
\end{equation}
as $t \to \infty$.

Conditioning on the event in \eqref{eqn:ndir-uniform} we have 

\begin{equation*}
\label{eqn:ndesc-bound}
n_A^{desc}(\cI_i) \preceq 
\sum_{j=1}^{ n_A^{dir} (\cI_i) } Y_j (\delta) \preceq 
\sum_{j=1}^{ C_1 W t \exp (t/3) } Y_j(\delta)
\end{equation*}
where $\{ Y_j(\cdot) : j \geq 1 \}$ is an sequence of $i.i.d.$ Yule process with rate $r^*_A$. Since $Y_j(t)$ follows a geometric distribution with parameter 
$\exp(-r^*_A t)$, we have that for any positive integer $C_2$:
\[
\pr (Y_j(\delta) \geq C_2) =
( 1-\exp(-r^*_A\delta) )^{C_2-1}
\leq 
(r^*_A \delta)^{C_2-1}
\]
Using a union bound it follows (for simplicity the conditional is suppressed) that
\begin{equation*}
\label{eqn:ndesc-tail}
\pr (n_A^{desc}(\cI_i)\geq C_1 C_2 W t e^{t/3}) \leq 
C_1 W t e^{t/3} \pr (Y_j(\delta) \geq C_2) \leq
C_1 W (r^*_A)^{C_2-1} t (e^{t/3})^{2-C_2}.
\end{equation*}
Apply union bound once again we get 
\[
\pr \left( \bigcup^{a / \delta}_{i=1} \Big \{ n_A^{desc}(\cI_i) \geq
C_1 C_2 W t e^{t/3} \Big \} \right) \leq
\frac{a}{\delta} \cdot C_1 W (r^*_A)^{C_2-1} t (e^{t/3})^{2-C_2} 
\to 0
\]
as $t \to \infty$ for $C_2 \geq 4$. Recall that $n_A^{bad}(\cI_i) \leq n_A^{desc}(\cI_i)$ and we can get rid of the conditional that was suppressed with the aid of \eqref{eqn:ndir-uniform}. Then letting $M = 4C_1$ completes the proof.

\qed

Now that we get the bad individuals under control, we can turn our attention to those good individuals. To start, combine Lemma \ref{lem:bad} with \eqref{eqn:na-uniform} we have
\begin{equation}
\label{eqn:good-uniform}
\pr \left( \bigcap^{ a / \delta }_{i=1} \Big\{ | n_A^{good}(\cI_i) - \delta W e^{t_i}| < 
4 \sqrt{t e^t} \Big\} \right) \to 1
\end{equation}
as $t \to \infty$, where $n_A^{good}(\cI_i)$ is the number of good individuals in $\cI_i$.

Next, note that good individuals in $\cI_i$ reproduce independently at rate $r^*_A$ starting from time $t_i + \delta = t_{i+1}$. In particular the probability that a good individual has $k$ offsprings by time $t$ is 
\[
g_k( t - t_{i+1} ) := \pr ( \text{Poisson}( r^*_A ( t - t_{i+1} ) ) = k).
\] 
Here $t - t_{i+1}$ refers to the time left until time $t$.
Since reproductions are all independent, we have (conditioning on $\BP(t_{i+1})$)
\begin{equation}
\label{eqn:good-deg-dist}
n_{k,A}^{good}(\cI_i) \equald \mathrm{Bin}(n_A^{good}(\cI_i), g_k( t - t_{i+1} ))
\end{equation}	
where $n_{k,A}^{good}(\cI_i)$ is the number of good individuals in $\cI_i$ that have $k$ offsprings by time $t$. Similar to the proof of \eqref{eqn:ndir-uniform}, by \eqref{eqn:good-uniform}, \eqref{eqn:good-deg-dist} and Hoeffding inequality for binomial distribution together with a union bound, there exists $C>0$ such that 
\begin{equation}
\label{eqn:good-deg-uniform}
\pr \left( \bigcap^{ a / \delta }_{i=1} \Big\{ |n_{k,A}^{good}(\cI_i) - 
\delta W e^{t_i} g_k( t - t_{i+1} )| < 
C W \log(W) t e^{t/3} \Big\} \right) 
\to 1
\end{equation}
as $t \to \infty$. 

With both good and bad individuals under control, we are ready to prove what we started out for. Note that
\begin{equation*}
\label{eqn:sandwich}
\sum_{i=1}^{ a / \delta } n_{k,A}^{good}(\cI_i) \leq 
n_{k,A} [t-a, t] \leq 
\sum_{i=1}^{ a / \delta } [ n_{k,A}^{good}(\cI_i) + n_A^{bad}(\cI_i) ].
\end{equation*}
For the bad individuals we know from Lemma \ref{lem:bad} that
\[
e^{-t} \sum_{i=1}^{ a / \delta } n_A^{bad}(\cI_i) \convp 0
\qquad 
\text{as}\qquad t \to \infty.
\]
On the other hand, for the good individuals we have from \eqref{eqn:good-deg-uniform} that
\begin{equation*}
\label{eqn:good-deg-integral}
e^{-t} \sum_{i=1}^{ a / \delta } n_{k,A}^{good}(\cI_i) - 
\sum_{i=1}^{ a / \delta } \delta W e^{-(t-t_i)} g_k( t - t_{i+1} ) 
\convp 0
\qquad 
\text{as}\qquad t \to \infty.
\end{equation*}
Finally, note that by definition of Riemann integral we know
\[
\sum_{i=1}^{ a / \delta } \delta W e^{-(t-t_i)} g_k( t - t_{i+1} ) \to
W \int_0^a \pr(\mathrm{Poisson}(r^*_A s)=k) e^{-s} ds 
\qquad 
\text{as} \qquad t \to \infty.
\]
This completes the proof for Proposition \ref{prop:trunc-deg-dist}. 

\qed

\subsection{Estimation of model parameters}
\label{sec:proof-est}
In this section we derive consistent estimators for model parameters $p$ and $q$ based on the limiting degree distribution (i.e. Corollary \ref{cor:estimators}).

Specifically, from Theorem \ref{thm:degree-dist} we know
\[
\frac{N_0(n)}{n} \convp \frac{p}{ 1 + r^*_A } + \frac{ 1-p }{ 1 + r^*_B }
\] 
and
\[
\frac{N_1(n)}{n} \convp 
\frac{p r^*_A}{ (1 + r^*_A)^2 } + 
\frac{ (1-p) r^*_B }{ (1 + r^*_B)^2 }.
\]
Let $\theta_A := 1 / ( 1 + r^*_A )$ and $\theta_B := 1 / ( 1 + r^*_B )$. If we can solve for unique $\hat{p}$, $\hat{\theta}_A$ and $\hat{\theta}_B$ (up to a switch between two types) that satisfy the following equations
\begin{equation}
\label{eqn:estimator-1}
m_1 := \frac{N_0(n)}{n} = \hat{p} \hat{\theta}_A + (1-\hat{p}) \hat{\theta}_B
\end{equation}

\begin{equation}
\label{eqn:estimator-2}
m_2 := \frac{ N_0(n) - N_1(n) }{n} = 
\hat{p} \hat{\theta}_A^2 + 
(1-\hat{p}) \hat{\theta}_B^2,
\end{equation}
then $\hat{p}$, $\hat{\theta}_A$ and $\hat{\theta}_B$ are consistent estimators of $p$, $\theta_A$ and $\theta_B$. We can further solve for a consistent estimator $\hat{q}$ for $q$ from either $\hat{\theta}_A$ or $\hat{\theta}_B$ together with $\hat{p}$. Note that all estimators are consistent because (as we shall see) they are continuous functions of $m_1$ and $m_2$.

To solve three unknowns from two equations \eqref{eqn:estimator-1} and \eqref{eqn:estimator-2}, we need a third equation
\begin{equation}
\label{eqn:estimator-3}
\frac{ \hat{p} }{ \hat{\theta}_A } + \frac{ 1-\hat{p} }{ \hat{\theta}_B } = 2 =
p ( 1 + r^*_A ) + (1-p) ( 1 + r^*_B )
\end{equation} 
where the second equality can be verified from the definitions of $r^*_A$ and $r^*_B$.

From \eqref{eqn:estimator-1} we have 
\begin{equation}
\label{eqn:estimator-p}
\hat{p} = \frac{ m_1 - \hat{\theta}_B }{ \hat{\theta}_A - \hat{\theta}_B }. 
\end{equation}
Plug it back to \eqref{eqn:estimator-2} and \eqref{eqn:estimator-3} we get 
$m_1 ( \hat{\theta}_A+ \hat{\theta}_B ) - \hat{\theta}_A \hat{\theta}_B = m_2$
and
$\hat{\theta}_A + \hat{\theta}_B - m_1 = 2 \hat{\theta}_A \hat{\theta}_B$.
Then we can solve for 
\[
\hat{\theta}_A \hat{\theta}_B = \frac { m_2 - m_1^2 }{ 2m_1 - 1 }
\quad \text{and} \quad
\hat{\theta}_A + \hat{\theta}_B  = \frac{ 2m_2 - m_1 }{ 2m_1 - 1 },
\]
so $\hat{\theta}_A$ and $\hat{\theta}_B$ are roots of the following quadratic equation:
\[
x^2 - \frac{ 2m_2 - m_1 }{ 2m_1 - 1 } x + \frac { m_2 - m_1^2 }{ 2m_1 - 1 } = 0.
\]
After we solve for $\hat{\theta}_A$ and $\hat{\theta}_B$ we can compute $\hat{p}$ from \eqref{eqn:estimator-p} and 
\[
\hat{q} = 
\frac{ \hat{p} - \hat{\theta}_A }{ (2\hat{p}-1) \hat{\theta}_A } = 
\frac{ 1 - \hat{p} - \hat{\theta}_B }{ (1-2\hat{p}) \hat{\theta}_B }
\]
from definitions of $\theta_A$ and $\theta_B$.

Note that these estimators fail when $m_1 - 1/2$, $m_2 - m_1^2$ or 
$2m_2 - m_1$ is negative. However, these cases are unlikely to occur for CMRT as all three quantities have positive limits as $n \to \infty$ (recall that $m_1$ and $m_2$ both depend on $n$). When $m_1 = 1/2$, these estimators also fail: from the remark under Theorem \ref{thm:degree-dist} we know $m_1$ has limit $1/2$ if and only if $p=1/2$ or $q=1$, and parameters $p$, $q$ are not identifiable from the limiting degree distribution in these special cases.

\subsection{Proof for maximal degree}
\label{sec:proof-max-deg}
In this section we shall prove Theorem \ref{thm:max-degree}. Throughout the proof we work with the continuous time embedding unless otherwise noted.

First we consider cases where $p,q \neq 1$. By Lemma \ref{lem:stop-time} and Egoroff's Theorem, given any $\epsilon>0$, we can choose $K>0$ such that (the dependence of $K$ on $\epsilon$ is suppressed throughout):
\begin{equation}
\label{eqn:stop-time-asy}
\pr \left( \sup_n |T_n - \log{n}| < K \right) > 1- \epsilon.
\end{equation}
To ease notations, write $T_n^+ = \log{n} + K$ and $T_n^- = \log{n} - K$.

{\bf Lower bound:} Here it is enough to consider just the type $A$ root. Recall from Section \ref{sec:ct-em} that this vertex reproduces at constant rate 
\begin{equation}
\label{eqn:r-def}
r^*_A = q + \frac{ (1-p)(1-q) }{p}.
\end{equation} 
Denote by $D(t)$ the out-degree of type $A$ root at time $t$. Then $D(t)$ follows a $\mbox{Poisson}(r^*_A t)$ distribution. So for any $0<\gamma<1$, by standard tail bound for Poisson distribution we have
\begin{equation}
\label{eqn:max-deg-lower-1}
\pr\left( D(T_n)\leq \gamma r^*_A \log{n} \right) \leq 
\pr\left( D(T_n^-)\leq \gamma r^*_A \log{n} \right) 
\leq \exp \{ - M \log{n} \} 
\end{equation} 
for large enough $n$, conditioning on the event in \eqref{eqn:stop-time-asy}. Here $M>0$ is a constant that depends on both $\gamma$ and $r^*_A$. 

Denote by $M(t)$ the maximal number of offsprings an individual has by time $t$. Since $M(t) \geq D(t)$, it follows from \eqref{eqn:max-deg-lower-1} that
\begin{equation}
\label{eqn:max-deg-lower-2}
\pr\left( M(T_n)\leq \gamma r^*_A \log{n} \right) 
\leq n^{-M}
\end{equation} 
conditioning on the event in \eqref{eqn:stop-time-asy}.

Let $\{ n_k \}_{k \geq 1}$ be an increasing sequence of positive integers such that $n_k \geq k^{2/M}$. Then from \eqref{eqn:max-deg-lower-2} we know
\begin{equation}
\label{eqn:max-deg-lower-3}
\pr\left( M(T_{n_k})\leq \gamma r^*_A \log{n_k} \right) 
\leq \frac{1}{k^2}
\end{equation} 
for large $k$, conditioning on the event in \eqref{eqn:stop-time-asy}. By Borel-Cantelli lemma and \eqref{eqn:stop-time-asy} it follows that 
\[ 
\pr \left( \liminf_{k \to \infty}  \frac{ M(T_{n_k}) }{\log n_k} \leq 
\gamma r^*_A \right) \leq \eps.
\]
From that we can also get
\[ 
\pr \left( \liminf_{n \to \infty}  \frac{ M(T_n) }{\log n} \leq 
\gamma r^*_A \right) \leq \eps.
\]

Recall that $M_n \stackrel{d}{=} M(T_n)$ from the embedding, let $\gamma \to 1$ and $\eps \to 0$ we have
\[ 
\pr \left( \liminf_{n \to \infty}  \frac{ M_n }{\log n} \geq 
r^*_A \right) = 1.
\]

This completes the proof for the lower bound part, with $C_1 = r^*_A$.

Alternatively, we provide here another proof for the lower bound part that works when $q \neq 0$. Ideas used in the proof shall become useful later on.

Consider type $A$ individuals alone. We call a type $A$ individual ``pure-blooded'' if all of its antecedents are type $A$ individuals. Define $\BP_A(t)$ to be the branching process consisting of all pure-blooded type $A$ individuals in $\BP(t)$. Let $T_{n,A} := \inf\set{ t \geq 0 : |\BP_A(t)| =n }$. Note that pure-blooded type $A$ individuals give birth to new pure-blooded type $A$ individuals at constant rate $q$, mimicking the results for CMRT (\ie Lemma \ref{lem:embed} and Lemma \ref{lem:stop-time} ) it is not hard to see that:
\begin{itemize}
	\item $\set{\BP_A(T_{n,A})}_{n\geq 1}  \stackrel{d}{=} \set{\cU_n}_{n\geq 1}$ as processes, where $\set{\cU_n}_{n\geq 1}$ is a URT.  
	\item The sequence of stopping times $T_{n,A}$ satisfy
	\begin{equation}
	\label{eqn:pureA-stop-time}
	T_{n,A} - \frac{1}{q}\log{n} \convas -\log(W_A)
	\end{equation}
	for a finite positive random variable $W_A\stackrel{d}{=} \exp(1)$. 
\end{itemize} 
Since we can transform results for URT into that of the branching process $\BP_A(t)$ through embedding, using Devroye and Lu \cite{devroye1995strong} we get
\[
\frac{M_A^{pure} (T_{n,A})}{\log n} \convas \frac{1}{\log 2}
\]
where $M_A^{pure}(t)$ is the maximal number of offsprings a pure-blooded type $A$ individual has by time $t$. Similar to \eqref{eqn:stop-time-asy} we have that by \eqref{eqn:pureA-stop-time} and Egoroff's Theorem, given any $\epsilon>0$, we can choose $K_A>0$ such that (again for simplicity the dependence of $K_A$ on $\epsilon$ is suppressed throughout):
\begin{equation}
\label{eqn:pureA-stop-time-asy}
\pr \left( \sup_n  |T_{n,A} -  \frac{1}{q}\log{n}| < K_A \right) > 1 - \epsilon.
\end{equation}
For each $t>0$, let $n^*(t)$ be the positive integer satisfying 
$ \frac{1}{q} \log(n^*(t)) + K_A \leq t <  \frac{1}{q} \log(n^*(t)+1) + K_A$. Conditioning on the event in \eqref{eqn:pureA-stop-time-asy} we have
\[
\frac{ M_A^{pure}(t) }{t} \geq 
\frac{ M_A^{pure}( \frac{1}{q} \log(n^*(t)) + K_A) }{  \frac{1}{q} \log(n^*(t)+1) + K_A }  \geq
\frac{ M_A^{pure}(T_{n^*(t),A}) }{  \frac{1}{q} \log(n^*(t)+1) + K_A }
\convas 
\frac{q}{\log 2}
\]
as $t \to \infty$, and it follows that 
\[
\pr \left( \liminf_{t \to \infty}  \frac{ M_A^{pure}(t) }{t} \geq 
\frac{q}{\log 2} \right) > 
1 - \epsilon.
\]
Denote by $M(t)$ the maximal number of offsprings an individual has by time $t$. Let $\epsilon \to 0$ in the above inequality and note $M(t) \geq M_A^{pure}(t)$, we have
\[ 
\pr \left( \liminf_{t \to \infty}  \frac{ M(t) }{t} \geq \frac{q}{\log 2} \right) = 1.
\]
Finally, since 
\begin{equation}
\label{eqn:stop-time-scale}
\frac{T_n}{\log n} \convas 1,
\end{equation}
by Lemma \ref{lem:stop-time}, and $M_n \stackrel{d}{=} M(T_n)$ from the embedding, we get
\[ 
\pr \left( \liminf_{n \to \infty}  \frac{ M_n }{\log n} \geq 
\frac{q}{\log 2} \right) = 1.
\]
This completes the proof for the lower bound part when $q \neq 0$, with $C_1 = \frac{q}{\log 2}$.

{\bf Upper bound:} 
First we consider the simpler type $A$ individuals. Given 
$\set{\BP(t)}_{t\geq 0}$, we couple it with another process where: whenever a type $A$ individual is born to a type $B$ individual in $\set{\BP(t)}_{t\geq 0}$, chose uniformly at random a living type $A$ individual, and treat the newborn individual as an offspring of that chosen individual. In this new process, if we look at type $A$ individuals alone, it is not hard to see that they give birth to new type $A$ individuals at constant rate 1. Denote by $M_A(t)$ and $\widetilde{M}_A(t)$ respectively the maximal number of offsprings a type $A$ individual has by time $t$ in $\set{\BP(t)}_{t\geq 0}$ and new process. Using the same argument as in the proof for lower bound we see that
\[ 
\pr \left( \limsup_{t \to \infty}  \frac{ \widetilde{M}_A(t) }{t} \leq 
\frac{1}{\log 2} \right) = 1.
\]
Since $\widetilde{M}_A(t) \geq M_A(t)$, this immediately implies
\[ 
\pr \left( \limsup_{t \to \infty}  \frac{ M_A(t) }{t} \leq 
\frac{1}{\log 2} \right) = 1.
\]
Denote by $M_{n,A}$ the maximal degree of type $A$ vertices in $\cT_n$.
Then it follows from \eqref{eqn:stop-time-scale} and $M_{n,A} \stackrel{d}{=} M_A(T_n)$ that
\[ 
\pr \left( \limsup_{n \to \infty}  \frac{ M_{n,A} }{\log n} \leq 
\frac{1}{\log 2} \right) = 1.
\]

Next we consider the more complicated type $B$ individuals. For a given type $B$ individual $v$ and time $T\in [0, t)$, write $n_v [T,t]$ for the number of offsprings this individual produced in time interval $[T,t]$. This is a pure birth process with rate 
\begin{equation*}
\label{eqn:c-def}
c \frac{n_A(t)} {n_B(t)} : =
\left( 1-q + \frac{ (1-p) q }{p} \right) \frac{n_A(t)}{n_B(t)}. 
\end{equation*}
Therefore for a fixed $T$, the following process 
\begin{equation}
\label{eqn:xt-def}
X(t) = n_v [T,t] - \int_T^t c \frac{n_A(s)}{n_B(s)}ds, \qquad t \geq T,
\end{equation}
is a martingale (here the infinitesimal generator is exactly the rate). 

To handle the variability of $X(t)$, we will need its predictable quadratic variation process $\langle X \rangle (t)$. Note that 
\[
\cA X^2(t) = \cA n_v^2 [T,t] + 
\cA \left( \int_T^t c \frac{n_A(s)}{n_B(s)} ds \right)^2 - 
2 \cA \left( n_v [T,t] \int_T^t c \frac{n_A(s)}{n_B(s)} ds \right)
\] 
with 
\[
\cA n_v^2 [T,t] = 
(2 n_v [T,t] + 1) c \frac{n_A(t)}{n_B(t)},
\]
\[
\cA \left( \int_T^t c \frac{n_A(s)}{n_B(s)} ds \right)^2 = 
2 (\int_T^t c \frac{n_A(s)}{n_B(s)} ds) c \frac{n_A(t)}{n_B(t)}
\]
and 
\[
\cA \left( n_v [T,t] \int_T^t c \frac{n_A(s)}{n_B(s)} ds \right) = 
n_v [T,t] c \frac{n_A(t)}{n_B(t)} + 
(\int_T^t c \frac{n_A(s)}{n_B(s)} ds) c \frac{n_A(t)}{n_B(t)}.
\]

It then follows from some elementary algebra that $\cA X^2(t) = c n_A(t) / n_B(t)$ and
\begin{equation}
\label{eqn:xx-def}
\langle X \rangle (t) = \int_T^t c \frac{n_A(s)}{n_B(s)} ds, 
\qquad t \geq T.  
\end{equation}
Now use Lemma \ref{lem:difference} to choose $T$ such that 
\begin{equation}
\label{eqn:T-def}
\pr\left(\sup_{t\geq T}\left|\frac{n_A(t)}{n_B(t)} - \frac{p}{1-p}\right| > \eps\right) \leq \eps.
\end{equation}
Also, define the stopping time
\[
S = \inf \set{t \geq T : \left| \frac{n_A(t)}{n_B(t)} - \frac{p}{1-p} \right| > \eps}.
\]
Observe that by our choice of $T$ we have $\pr ( S < \infty ) \leq \eps$. The idea here is to bound $n_A(t) / n_B(t)$ around $p / (1-p)$ after some finite time $T$, and show that what happened before time $T$ does not have a noticeable effect in the long run.  

Recall that we write $T_n^+ = \log{n} + K$ and $T_n^- = \log{n} - K$. Consider the process $\set{X (t \wedge S): t \geq T}$ and note that for $n$ large enough we have $T_n^- \geq T$. By the exponential martingale inequality from Liptser and Shiryayev \cite[Section 4.13, Theorem 5]{liptser2012theory} with choices
\[
K=2 
\quad \text{and} \quad
\varphi(t) = c \left( \frac{p}{1-p} +\eps \right) t,
\]
we have for any $\delta>0$,
\begin{align}
\label{eqn:tail-bound} 
\pr \left( \sup_{t \leq T_n^+ \wedge S} X(t) \geq 
\delta c \left( \frac{p}{1-p} + \eps \right) T_n^+ \right) & \leq 
\exp \left( -\kappa c \left( \frac{p}{1-p} + \eps \right) T_n^+ \right) \notag \\
& + \pr \left( \langle X \rangle (T_n^+ \wedge S) \geq 2 \varphi (T_n^+) \right), 
\end{align}
where $\kappa = (\delta+2) \log{ \frac{\delta+2}{2} } - \delta$. By definition of $S$ and the expression of $\langle X \rangle (\cdot)$ in \eqref{eqn:xx-def}, we have that with probability one 
\[
\langle X \rangle (T_n^+ \wedge S) \leq \varphi (T_n^+),
\]
so the second term on the right hand side of \eqref{eqn:tail-bound} vanishes.
Further, using the expression of $X(\cdot)$ in \eqref{eqn:xt-def} we get
\begin{equation}
\label{eqn:tail-bound-h-stop}
\pr \left( n_v [T, T_n^+ \wedge S] \geq 
(\delta+1) c \left( \frac{p}{1-p} + \eps \right) T_n^+ \right) \leq
\exp \left( -\kappa c \left( \frac{p}{1-p} + \eps \right) T_n^+ \right). 
\end{equation}

For $a<b$, denote by $M_B [a, b]$ the maximal number of offsprings produced by a type $B$ individual in time interval $[a,b]$ and let 
$M_B(t) = M_B [0, t]$. Then using a union bound and note that 
$M_B [T, T_n \wedge S] \leq M_B [T, T_n^+ \wedge S]$ with probability at least $1-\eps$, we have 
\begin{equation*}
\label{eqn:tail-bound-mb-stop}
\pr \left( M_B [T, T_n \wedge S] \geq 
(\delta+1) c \left( \frac{p}{1-p} + \eps \right) T_n^+ \right) \leq
\eps + n \exp \left( -\kappa c \left( \frac{p}{1-p} + \eps \right) T_n^+ \right). 
\end{equation*}
Further, from our choice of $T$ in \eqref{eqn:T-def} it follows that
\begin{equation*}
\label{eqn:tail-bound-mb-1}
\pr \left( M_B [T, T_n] \geq 
(\delta+1) c \left( \frac{p}{1-p} + \eps \right) T_n^+ \right) \leq
2\eps + n \exp \left( -\kappa c \left( \frac{p}{1-p} + \eps \right) T_n^+ \right). 
\end{equation*}

Next, note that there exists $L>0$ such that
\begin{equation*}
\label{eqn:tail-bound-mb-2}
\pr \left( M_B(T) \geq L \right) < \eps. 
\end{equation*}
As $M_B(T_n) \leq M_B(T)+ M_B [T, T_n]$, combining the above results readily yields
\begin{equation}
\label{eqn:tail-bound-mb-3}
\pr \left( M_B (T_n) \geq 
L + (\delta+1) c \left( \frac{p}{1-p} + \eps \right) T_n^+ \right) \leq
3\eps + n \exp \left( -\kappa c \left( \frac{p}{1-p} + \eps \right) T_n^+ \right). 
\end{equation}
Denote by $M_{n,B}$ the maximal degree of type $B$ vertices in $\cT_n$. Letting $\eps \to 0$ in \eqref{eqn:tail-bound-mb-3} and note $M_{n,B} \stackrel{d}{=} M_B (T_n)$, we have
\begin{equation}
\label{eqn:tail-bound-mb-final}
\pr \left( M_{n,B} \geq 
L + (\delta+1) \frac{cp}{1-p} T_n^+ \right) \leq
n \exp \left( -\frac{\kappa c p}{1-p} T_n^+ \right). 
\end{equation}
Recall that $\kappa = (\delta+2) \log{ \frac{\delta+2}{2} } - \delta$. So given any $\widetilde{\eps} >0$, for 
$\delta > \max \{ 2e^2-2, \frac{ (2+\widetilde{\eps}) (1-p) }{ cp } - 4\}$ 
we have $\frac{\kappa c p}{1-p} T_n^+ > (2+\widetilde{\eps}) \log{n}$ for large enough $n$, which makes the right hand side of \eqref{eqn:tail-bound-mb-final} summable. By Borel-Cantelli lemma this implies 
\[
\pr \left( \limsup_{n \to \infty}  \frac{ M_{n,B} }{\log n} \leq 
(\delta+1) \frac{cp}{1-p}  + \widetilde{\eps} \right) = 1.
\]
Letting $\widetilde{\eps} \to 0$ we get 
\[
\pr \left( \limsup_{n \to \infty}  \frac{ M_{n,B} }{\log n} \leq C_2 \right) = 1
\]
where $C_2 = \max \{ \frac{ (2e^2-1)cp }{1-p}, 2 - \frac{3cp}{1-p} \}$.
This completes the proof for the upper bound part.

This left us with the cases where $p=1$ or $q=1$. When $q=1$, as noted in Section \ref{sec:special}, $\set{\cT_n}_{n\geq 2}$ looks like two disjoint URTs connected by a single edge between roots. By strong law of large numbers, with probability one these two subtrees have sizes proportional to $p$ and $1-p$ asymptotically. Therefore, appealing to existing results on maximal degree of URT \cite{devroye1995strong} we have
\[
\frac{ M_n }{\log n} \convas \frac{1}{\log 2}.
\]

Last we consider the special case where $p = 1$ but $q \neq 1$. We shall use the same notations as defined before. First, note that what we proved for type $A$ individuals still holds. Specifically, recall that
\[ 
\pr \left( \limsup_{n \to \infty}  \frac{ M_{n,A} }{\log n} \leq 
\frac{1}{\log 2} \right) = 1.
\]
However, we can no longer define $T$ by \eqref{eqn:T-def}. In fact, consider the unique type $B$ vertex (\ie the type $B$ root) in $\cT_n$ and note that all type $A$ vertices other than the root have a fixed probability $1-q$ to connect to the type $B$ root independently. By strong law of large numbers applied to binomial distribution we have
\[ 
\frac{ M_{n,B} }{n} \convas 1-q,
\quad \text{and therefore} \quad
\frac{ M_n }{n} \convas 1-q.
\]

This completes the proof for Theorem \ref{thm:max-degree}.

\qed

\subsection{Proof for height}
\label{sec:proof-height}
In this section we shall prove Theorem \ref{thm:height}. Instead of proving the result from scratch using continuous time martingales as what we did for maximal degree, we present here a proof that appeals to existing results on first birth problem of branching processes. 

Once again, we consider first cases where $p,q \neq 1$. The basic idea here is still the same, \ie to bound $n_A(t) / n_B(t)$ around $p / (1-p)$ after some finite time and prove that what happened in the beginning does not really matter in the long term.

To obtain strong convergence, it is enough to prove that for any 
$\delta, \eps>0$,  
\begin{equation}
\label{eqn:h-str-conv}
\limsup_{n \to \infty} \pr \left( \sup_{k \geq n} 
|\frac{H_k}{\log k} - e| > \delta \right) < \eps
\end{equation}
Define the event 
\begin{equation}
\label{eqn:good-event-1}
E_1 = \Bigg\{ \sup_{t \geq T} \left| \frac{n_A(t)}{n_B(t)} - \frac{p}{1-p} \right| <
\eta \frac{p}{1-p} \Bigg\}.
\end{equation}
where $0<\eta<1$ is any given constant. By Lemma \ref{lem:difference} there exists $T>0$ such that $\pr(E_1) > 1 - \eps / 3$. Also, choose $N \in \bZ^+$ such that $\pr \left( n(T) > N \right) < \eps / 3$, and define another event $E_2 = \{ n(T) \leq N \}$. Moreover, choose $T_n^+$ and $T_n^-$ in the same way as in Section \ref{sec:proof-max-deg}, with $\eps$ replaced by $\eps / 3$, and define our last ``good'' event $E_3 = \{ T_n^- < T_n < T_n^+ \text{ for all n}\in \bZ^+ \}$. In what follows, we condition on the event $E_1 \cap E_2 \cap E_3$ and note that $\pr ( E_1 \cap E_2 \cap E_3 ) > 1 - \eps$ (for simplicity the conditional is suppressed throughout). 

On this event, we have at most $N$ individuals alive at time $T$ and the ratio of $n_A(t)$ to $n_B(t)$ is bounded around $p / (1-p)$ after that time. For a fixed individual $v$ alive at time $T$ and $t>0$, denote by $H_v(t)$ the height of the subtree root at $v$ in $\BP(T+t)$. To bound $H_v(t)$, we now construct two processes. Recall our continuous time process as defined in Section \ref{sec:ct-em}, and consider a process where each type $B$ individual gives birth to type $A$ individuals at rate $(1-\eta)p(1-q) / (1-p)$ and type $B$ individuals at rate $(1-\eta) q$, while everything else stays the same. Denote by $H_{min}(t)$ the height of this tree at time $t$, $B_{min}(n)$ the time when the first individual in the $n$-th generation is born, and define $H_{max}(t)$ similarly using $1+\eta$ instead of $1-\eta$ in the rates. Since the reproduction rates of type $A$ individuals are constants, and those of type $B$ individuals only depend on 
$n_A(t) / n_B(t)$, we have $H_{min}(t) \preceq H_v(t) \preceq H_{max}(t)$ where $\preceq$ denotes stochastic dominance as processes (\ie for any monotone increasing functional $f$ we have $f(H_{min}(t)) \preceq f(H_v(t)) \preceq f(H_{max}(t))$ where $\preceq$ denotes the usual stochastic dominance).

From Biggins \cite[Theorem 2]{biggins1976first} we know that 
\[
\lim_{n \to \infty} \frac{B_{min}(n)}{n} = \gamma_{min},
\] 
where $\gamma_{min}$ can be calculated following the procedure given in the paper. First, compute the matrix $\Phi(\theta)$ with entries 
\[
\Phi_{ij} (\theta) = \theta \int_{0}^{\infty} e^{-\theta t} \bE (Z_{ij}(t)) dt = 
\frac{r_{ij}}{\theta}.
\] 
Here $Z_{ij}(t)$ denotes the number of type $j$ individuals born to a type $i$ individual by time $t$, and $r_{ij}$ denotes the rate at which a type $i$ individual gives birth to type $j$ individuals. Then take the largest eigenvalue 
$\phi (\theta) = \lambda_{min} / \theta$ of $\Phi(\theta)$. In our case we have
\[
\lambda_{min} = \frac{ (2-\eta)q + \sqrt{ \eta^2 q^2 + 4(1-\eta)(1-q)^2 } }{2}
\]
by calculation. Finally, define 
\[
\mu(a) = \inf \Bigg\{ e^{\theta a} \phi(\theta): \theta > 0 \Bigg\},
\] 
and compute
$\gamma_{min} = \inf \Big\{ a: \mu(a) \geq 1 \Big\} = 1 / (\lambda_{min} e)$.

Since
\[
B_{min} (H_{min}(t)) \leq t \leq B_{min} (H_{min}(t)+1),
\] 
dividing by $H_{min}(t)$ and letting $t \to \infty$ we get
\begin{equation}
\label{eqn:h-min-as}
\frac{H_{min}(t)}{t} \convas \lambda_{min} e.
\end{equation} 

Similarly we have
\begin{equation}
\label{eqn:h-max-as}
\frac{H_{max}(t)}{t} \convas \lambda_{max} e,
\end{equation} 
where 
\[
\lambda_{max} = \frac{ (2+\eta)q + \sqrt{ \eta^2 q^2 + 4(1+\eta)(1-q)^2 } }{2}.
\]
As the eigenvalues are continuous with respect to $\eta$, and 
$\lambda_{max} = \lambda_{min} = 1$ when $\eta=0$, we can choose $\eta$ in \eqref{eqn:good-event-1} small enough such that both $\lambda_{max} - 1$ and $1 - \lambda_{min}$ are smaller than $\delta / 3$.

\begin{rem}
	Since we will need to generalize the result to CMRT with more types, we include here an alternative argument using Perron-Frobenius theory of positive matrices \cite[Chapter 8]{meyer2000matrix}. Note that when $\eta=0$, the rate matrix consisting of $r_{ij}$'s is a positive matrix (\ie a matrix where all entries are strictly positive)
	\[
	\left( \begin{array}{cc}
	q  &  \frac{(1-p)(1-q)}{p}\\
	\frac{p(1-q)}{1-p} & q \\
	\end{array} \right)
	\]
	with left eigenvector $(p , 1-p)$ corresponding to eigenvalue $\lambda = 1$. Since this eigenvector has strictly positive coordinates, we have that $\lambda = 1$ is the unique largest eigenvalue of that rate matrix. Therefore the largest eigenvalue of rate matrix is continuous with respect to $\eta$ around $\eta=0$, and the previous result we established on $\lambda_{max}$ and $\lambda_{min}$ follows.
\end{rem}

By \eqref{eqn:h-min-as}, \eqref{eqn:h-max-as} and stochastic dominance (as processes) we see that 
\[
\lim_{s \to \infty} \pr \left( \inf_{t \geq s} \frac{H_{v}(t)}{t} < 
\lambda_{min} e - \frac{\delta}{3} \right) = 0 
\quad \text{and} \quad
\lim_{s \to \infty} \pr \left( \sup_{t \geq s} \frac{H_{v}(t)}{t} > 
\lambda_{max} e + \frac{\delta}{3} \right)=0
\] 
hold for all $v \in V$, where $V$ denotes the set of individuals alive at time $T$. 

Since $H_n \geq H_v (T_n^- -T)$, we have 
\[
\limsup_{n \to \infty} \pr \left( \inf_{k \geq n} 
\frac{H_k}{T_k^- -T} < \lambda_{min} e - \frac{\delta}{3} \right) \leq 
\lim_{s \to \infty} \pr \left( \inf_{t \geq s} \frac{H_{v}(t)}{t} < 
\lambda_{min} e - \frac{\delta}{3} \right) = 0.
\] 

On the other hand, note that $H_n \leq \max_{u\in V} H_u (T_n^+ -T) + N$. Using union bound with $|V| \leq N$ we get
\[
\limsup_{n \to \infty} \pr \left( \sup_{k \geq n} 
\frac{H_k - N}{T_k^+ -T} > \lambda_{max} e + \frac{\delta}{3} \right) \leq 
\sum_{u \in V} \lim_{s \to \infty} \pr \left( 
\sup_{t \geq s} \frac{H_{u}(t)}{t} > \lambda_{max} e + \frac{\delta}{3} \right) 
=0.
\]
With
\[
| \frac{H_k}{T_k^- -T} - \frac{H_k}{\log k} | < \frac{\delta}{3}
\quad \text{and} \quad
| \frac{H_k - N}{T_k^+ -T} - \frac{H_k}{\log k} | < \frac{\delta}{3}
\] 
for large enough $k$, by our choice of $\eta$ and triangle inequality it follows that
\begin{equation}
\label{eqn:h-str-conv-cd}
\limsup_{n \to \infty} \pr \left( \sup_{k \geq n}
| \frac{H_k}{\log k} - e |> \delta \right) = 0.
\end{equation}
Finally, do not forget that we are conditioning on the event 
$E_1 \cap E_2 \cap E_3$, which occurs with probability at least $1 - \eps$. Therefore we have \eqref{eqn:h-str-conv} as desired.

Once again we are left with cases where $p=1$ or $q=1$. When $q=1$, recall that $\set{\cT_n}_{n\geq 2}$ looks like two disjoint URTs in this case. Using results for URT \cite{pittel1994note} we have
\[
\frac{H_n}{\log n} \convas e.
\]

Last we turn to cases where $p = 1$ but $q \neq 1$. Consider only pure-blooded type $A$ individuals as in Section \ref{sec:proof-max-deg} and use the same argument there we have 
\[ 
\pr \left( \liminf_{n \to \infty} \frac{H_n}{\log n} \geq qe \right) = 1.
\]
For upper bound we construct $\set{\cT_n}_{n\geq 2}$ from a URT $\set{\cU_n}_{n\geq 1}$ as described in special cases of Section \ref{sec:special} and note that by construction the height of $\cT_n$ is at most equal to that of $\cU_n$. From results for URT \cite{pittel1994note} we know that
\[ 
\pr \left( \limsup_{n \to \infty} \frac{H_n}{\log n} \leq e \right) = 1.
\]

This completes the proof for Theorem \ref{thm:height}.

\qed

\subsection{Extension to general CMRT}
\label{sec:proof-cmrt-general}

For CMRT with $K$ types and arbitrary attachment probabilities, we will need a population-dependent branching process defined as follows:
\begin{enumeratea}
	\item {\bf Initialization:} start with $K$ individuals at $t=0$, with one of each type. For any time $t\geq 0$ and $1 \leq i \leq K$, let $n_i(t)$ be the number of type $i$ individuals. We have $n_i(0) = 1$. Denote by $\cF(t)$ the $\sigma$-field generated by the process until time $t$. 
	\item {\bf Types:} Each individual in the system has a type $\in \set{1,2,...,K}$ and lives forever, while giving birth to individuals of all types.
	\item {\bf Reproduction:} At any time $t$, a living type $i$ individual gives birth to type $j$ ($1 \leq j \leq K$) individuals at rate: 
		\[
		r_{ij}(t) = \frac{n_1(t)}{n_i(t)} \cdot \frac{p_{j}q_{ji}}{p_1}.
		\]
\end{enumeratea}

Then using the same arguments as in the proofs for CMRT with two types we can prove similar results in the general case as stated in Section \ref{sec:cmrt-res-general}. 

\subsection{Extension to CWRT}
\label{sec:proof-cwrt}

Similar to CMRT, we can define a population-dependent branching process with $K$ types for CWRT, with birth rates replaced by
\[
	r_{ij}(t) = \frac{n_1(t)}{p_1} \cdot \frac{p_j \omega_{ji}}{\sum_{l=1}^K \omega_{jl} n_l(t)}.
\]

Then following the same procedure as in the proof for CMRT we can derive the limiting degree distribution of CWRT, \ie Theorem \ref{thm:degree-dist-cwrt}.

\subsection{Extension to CMPA}
\label{sec:proof-cmpa}

For CMPA the birth rates have a more complex form and we need to introduce some notations. 

For any time $t\geq 0$, $1 \leq i,j \leq K$ and individual $v$, let $n_{ij}(t)$ be the number of type $j$ individuals born to type $i$ individuals, and $d_v(t)$ be the number of offsprings born to individual $v$. Then at any time $t$, a living type $i$ individual $v$ gives birth to type $j$ ($1 \leq j \leq K$) individuals at rate: 
\begin{equation}
\label{eqn:rate-cmpa}
r_v(t,j) = \frac{n_1(t)}{p_1} \cdot \frac{p_{j}q_{ji} (d_v(t) + \alpha_{ji})}{\alpha_{ji} n_i(t) + \sum_{l=1}^K n_{il}(t)}.
\end{equation}
Note here that $\sum_{l=1}^K n_{il}(t)$ is exactly the total number of offsprings born to type $i$ individuals.

To derive the limiting degree distribution of CMPA, \ie Theorem \ref{thm:degree-dist-cmpa}, we can still mimic the proof for CMRT. However, computation for the integral in the final step will be different. Recall that for CMRT we have the integral 
\[
\int_0^\infty \pr(\mathrm{Poisson}(r^*_A s)=k) e^{-s} ds  = \frac{1}{ 1 + r^*_A } \left( \frac{r^*_A}{ 1 + r^*_A } \right)^k.
\]
For CMPA, the (time-homogeneous) Poisson process is replaced by a time-inhomogeneous Poisson process. To get the birth rates, observe that similar to the proof of CMRT we have
$
n_i(t)/n(t) \convas p_i
$
and 
$
n_{ij}(t)/n(t) \convas p_j q_{ji}
$
as $t \to \infty$. Therefore from \eqref{eqn:rate-cmpa} the reproduction of a type $i$ individual $v$ is approximately 
\[
\sum_{j=1}^K \frac{p_j q_{ji} (d_v(t) + \alpha_{ji}) } {\alpha_{ji} p_i + \sum_{l=1}^K p_l q_{li} } = \nu_i (d_v(t) + \alpha_i ).
\]
This leads to a pure birth process with rate $\nu_i (m + \alpha_i -1)$, where $m$ is the current population size. We shall denote this process by $\{ \tilde{Y}(\nu_i, t) : t \geq 0 \}$ (with $\tilde{Y}(\nu_i, 0) = 1$). Then similar to CMRT we need to compute the integral $\int_0^\infty \pr( \tilde{Y}(\nu_i, s) =k+1) e^{-s} ds$.

Since transition probability function of $\{ \tilde{Y}(\nu_i, t) : t \geq 0 \}$ can be computed explicitly (see \eg \cite[Proposition 6.1]{ross2014introduction}), we have
\[
\pr(\tilde{Y}(\nu_i, s) =k+1) = \frac{\Gamma(k+\alpha_i)}{\Gamma(\alpha_i)\Gamma(k+1)} e^{-\alpha_i \nu_i s} (1- e^{-\nu_i s})^{k}.
\]
Plugging it into the integral we get
\[
\int_0^\infty \pr(Y(\nu_i, s) =k+1) e^{-s} ds = \frac{\Gamma(k+\alpha_i)}{\Gamma(\alpha_i)\Gamma(k+1)} \int_0^\infty e^{-\alpha_i \nu_i s} (1- e^{-\nu_i s})^{k} e^{-s} ds.
\]
Let $x=e^{-\nu_i s}$ and by change of variable we have
\[
\int_0^\infty e^{-\alpha_i \nu_i s} (1- e^{-\nu_i s})^{k} e^{-s} ds = \frac{1}{\nu_i} \int_0^1 x^{1/\nu_i + \alpha_i -1} (1-x)^k dx = \frac{1}{\nu_i} \frac{\Gamma(k+1) \Gamma(1/\nu_i + \alpha_i)}{\Gamma(k+1+1/\nu_i + \alpha_i)}.
\]
This leads to the constant in \eqref{const-cmpa}.

\section{Acknowledgments}
The research of SB and RF was supported in part by the NSF (DMS-1606839, DMS-1613072)  and the Army Research Office (W911NF-17-1-0010). The research of RF was supported in part by the NSF under Grant DMS-1638521 to the Statistical and Applied Mathematical Sciences Institute. The research of AN was supported in part by NSF grants DMS-1613072 and DMS-161326, and by NIH grant R01 HG009125-01.

\bibliographystyle{plain}
\bibliography{comm_trees}

\begin{bibdiv}
\begin{biblist}

\bib{aldous1991asymptotic}{article}{
      author={Aldous, David},
       title={Asymptotic fringe distributions for general families of random
  trees},
        date={1991},
     journal={The Annals of Applied Probability},
       pages={228\ndash 266},
}

\bib{athreya1968embedding}{article}{
      author={Athreya, Krishna~B},
      author={Karlin, Samuel},
       title={Embedding of urn schemes into continuous time markov branching
  processes and related limit theorems},
        date={1968},
     journal={The Annals of Mathematical Statistics},
      volume={39},
      number={6},
       pages={1801\ndash 1817},
}

\bib{barabasi1999emergence}{article}{
      author={Barab{\'a}si, Albert-L{\'a}szl{\'o}},
      author={Albert, R{\'e}ka},
       title={Emergence of scaling in random networks},
        date={1999},
     journal={science},
      volume={286},
      number={5439},
       pages={509\ndash 512},
}

\bib{bhamidi2015change}{article}{
      author={Bhamidi, Shankar},
      author={Jin, Jimmy},
      author={Nobel, Andrew},
       title={Change point detection in network models: Preferential attachment
  and long range dependence},
        date={2015},
     journal={arXiv preprint arXiv:1508.02043},
}

\bib{biggins1976first}{article}{
      author={Biggins, JD},
       title={The first-and last-birth problems for a multitype age-dependent
  branching process},
        date={1976},
     journal={Advances in Applied Probability},
      volume={8},
      number={3},
       pages={446\ndash 459},
}

\bib{bollobas2001degree}{article}{
      author={Bollob{\'a}s, B{\'{}}~ela},
      author={Riordan, Oliver},
      author={Spencer, Joel},
      author={Tusn{\'a}dy, G{\'a}bor},
       title={The degree sequence of a scale-free random graph process},
        date={2001},
     journal={Random Structures \& Algorithms},
      volume={18},
      number={3},
       pages={279\ndash 290},
}

\bib{bordenave2015non}{inproceedings}{
      author={Bordenave, Charles},
      author={Lelarge, Marc},
      author={Massouli{\'e}, Laurent},
       title={Non-backtracking spectrum of random graphs: community detection
  and non-regular ramanujan graphs},
organization={IEEE},
        date={2015},
   booktitle={2015 ieee 56th annual symposium on foundations of computer
  science},
       pages={1347\ndash 1357},
}

\bib{bubeck2017finding}{article}{
      author={Bubeck, S{\'e}bastien},
      author={Devroye, Luc},
      author={Lugosi, G{\'a}bor},
       title={Finding adam in random growing trees},
        date={2017},
     journal={Random Structures \& Algorithms},
      volume={50},
      number={2},
       pages={158\ndash 172},
}

\bib{crump1968general}{article}{
      author={Crump, Kenny~S},
      author={Mode, Charles~J},
       title={A general age-dependent branching process. i},
        date={1968},
     journal={Journal of Mathematical Analysis and Applications},
      volume={24},
      number={3},
       pages={494\ndash 508},
}

\bib{deijfen2016birds}{article}{
      author={Deijfen, Maria},
      author={Fitzner, Robert},
       title={Birds of a feather or opposites attract-effects in network
  modelling},
        date={2016},
     journal={arXiv preprint arXiv:1612.03127},
}

\bib{devroye1998branching}{incollection}{
      author={Devroye, Luc},
       title={Branching processes and their applications in the analysis of
  tree structures and tree algorithms},
        date={1998},
   booktitle={Probabilistic methods for algorithmic discrete mathematics},
   publisher={Springer},
       pages={249\ndash 314},
}

\bib{devroye2012depth}{article}{
      author={Devroye, Luc},
      author={Fawzi, Omar},
      author={Fraiman, Nicolas},
       title={Depth properties of scaled attachment random recursive trees},
        date={2012},
     journal={Random Structures \& Algorithms},
      volume={41},
      number={1},
       pages={66\ndash 98},
}

\bib{devroye1995strong}{article}{
      author={Devroye, Luc},
      author={Lu, Jiang},
       title={The strong convergence of maximal degrees in uniform random
  recursive trees and dags},
        date={1995},
     journal={Random Structures \& Algorithms},
      volume={7},
      number={1},
       pages={1\ndash 14},
}

\bib{d2007power}{article}{
      author={D'Souza, Raissa~M},
      author={Krapivsky, Paul~L},
      author={Moore, Cristopher},
       title={The power of choice in growing trees},
        date={2007},
     journal={The European Physical Journal B},
      volume={59},
      number={4},
       pages={535\ndash 543},
}

\bib{fortunato2010community}{article}{
      author={Fortunato, Santo},
       title={Community detection in graphs},
        date={2010},
     journal={Physics reports},
      volume={486},
      number={3-5},
       pages={75\ndash 174},
}

\bib{gastwirth1977probability}{article}{
      author={Gastwirth, Joseph~L},
       title={A probability model of a pyramid scheme},
        date={1977},
     journal={The American Statistician},
      volume={31},
      number={2},
       pages={79\ndash 82},
}

\bib{gastwirth1984two}{article}{
      author={Gastwirth, Joseph~L},
      author={Bhattacharya, PK},
       title={Two probability models of pyramid or chain letter schemes
  demonstrating that their promotional claims are unreliable},
        date={1984},
     journal={Operations Research},
      volume={32},
      number={3},
       pages={527\ndash 536},
}

\bib{jagers1996asymptotic}{article}{
      author={Jagers, Peter},
      author={Nerman, Olle},
       title={The asymptotic composition of supercritial, multitype branching
  populations},
        date={1996},
     journal={LECTURE NOTES IN MATHEMATICS-SPRINGER VERLAG-},
       pages={40\ndash 54},
}

\bib{janson2005asymptotic}{article}{
      author={Janson, Svante},
       title={Asymptotic degree distribution in random recursive trees},
        date={2005},
     journal={Random Structures \& Algorithms},
      volume={26},
      number={1-2},
       pages={69\ndash 83},
}

\bib{karatzas2012brownian}{book}{
      author={Karatzas, Ioannis},
      author={Shreve, Steven},
       title={Brownian motion and stochastic calculus},
   publisher={Springer Science \& Business Media},
        date={2012},
      volume={113},
}

\bib{karrer2011stochastic}{article}{
      author={Karrer, Brian},
      author={Newman, Mark~EJ},
       title={Stochastic blockmodels and community structure in networks},
        date={2011},
     journal={Physical review E},
      volume={83},
      number={1},
       pages={016107},
}

\bib{krzakala2013spectral}{article}{
      author={Krzakala, Florent},
      author={Moore, Cristopher},
      author={Mossel, Elchanan},
      author={Neeman, Joe},
      author={Sly, Allan},
      author={Zdeborov{\'a}, Lenka},
      author={Zhang, Pan},
       title={Spectral redemption in clustering sparse networks},
        date={2013},
     journal={Proceedings of the National Academy of Sciences},
      volume={110},
      number={52},
       pages={20935\ndash 20940},
}

\bib{liptser2012theory}{book}{
      author={Liptser, Robert},
      author={Shiryayev, Albert~Nikolaevich},
       title={Theory of {M}artingales},
   publisher={Springer Science \& Business Media},
        date={2012},
      volume={49},
}

\bib{mahmoud2010power}{article}{
      author={Mahmoud, Hosam~M},
       title={The power of choice in the construction of recursive trees},
        date={2010},
     journal={Methodology and Computing in Applied Probability},
      volume={12},
      number={4},
       pages={763\ndash 773},
}

\bib{mahmoud1992asymptotic}{article}{
      author={Mahmoud, Hosam~M},
      author={Smythe, Robert~T.},
       title={Asymptotic joint normality of outdegrees of nodes in random
  recursive trees},
        date={1992},
     journal={Random Structures \& Algorithms},
      volume={3},
      number={3},
       pages={255\ndash 266},
}

\bib{meyer2000matrix}{book}{
      author={Meyer, Carl~D},
       title={Matrix analysis and applied linear algebra},
   publisher={Siam},
        date={2000},
      volume={71},
}

\bib{moon1974distance}{article}{
      author={Moon, JW},
       title={The distance between nodes in recursive trees},
        date={1974},
     journal={London Mathematical Society Lecture Notes},
      volume={13},
       pages={125\ndash 132},
}

\bib{najock1982number}{article}{
      author={Najock, Dietmar},
      author={Heyde, CC},
       title={On the number of terminal vertices in certain random trees with
  an application to stemma construction in philology},
        date={1982},
     journal={Journal of Applied Probability},
      volume={19},
      number={3},
       pages={675\ndash 680},
}

\bib{newman2012communities}{article}{
      author={Newman, Mark~EJ},
       title={Communities, modules and large-scale structure in networks},
        date={2012},
     journal={Nature physics},
      volume={8},
      number={1},
       pages={25\ndash 31},
}

\bib{oksendal2003stochastic}{incollection}{
      author={{\O}ksendal, Bernt},
       title={Stochastic differential equations},
        date={2003},
   booktitle={Stochastic differential equations},
   publisher={Springer},
       pages={65\ndash 84},
}

\bib{pittel1994note}{article}{
      author={Pittel, Boris},
       title={Note on the heights of random recursive trees and random m-ary
  search trees},
        date={1994},
     journal={Random Structures \& Algorithms},
      volume={5},
      number={2},
       pages={337\ndash 347},
}

\bib{rosengren2017multi}{article}{
      author={Rosengren, Sebastian},
       title={A multi-type preferential attachment tree},
        date={2017},
     journal={arXiv preprint arXiv:1704.03256},
}

\bib{ross2014introduction}{book}{
      author={Ross, Sheldon~M},
       title={Introduction to probability models},
   publisher={Academic press},
        date={2014},
}

\bib{smythe1995survey}{article}{
      author={Smythe, Robert~T},
      author={Mahmoud, Hosam~M},
       title={A survey of recursive trees},
        date={1995},
     journal={Theory of Probability and Mathematical Statistics},
      number={51},
       pages={1\ndash 28},
}

\end{biblist}
\end{bibdiv}

\end{document}